\documentclass[11pt,reqno]{amsart}

\usepackage{amsmath,amssymb,mathrsfs}
\usepackage{graphicx,cite,times}

\numberwithin{equation}{section}

\newcommand{\al}{\alpha}
\newcommand{\ba}{\beta}
\newcommand{\ga}{\gamma}
\newcommand{\de}{\delta}
\newcommand{\om}{\omega}

\newcommand{\ve}{\varepsilon}

\newcommand{\te}{\theta}
\newcommand{\ka}{\kappa}
\newcommand{\la}{\lambda}
\newcommand{\rh}{\rho}

\newcommand{\ps}{\psi}
\newcommand{\ph}{\phi}

\newcommand{\nt}{\notag}
\newcommand{\lb}{\label}

\newcommand{\Ga}{\Gamma}
\newcommand{\De}{\Delta}

\newcommand{\La}{\Lambda}
\newcommand{\Om}{\Omega}

\newcommand{\tl}{\tilde}
\newcommand{\qu}{\quad}
\newcommand{\na}{\nabla}
\newcommand{\hs}{\hspace}

\newcommand{\mb}{\mathbb}
\newcommand{\mc}{\mathcal}
\newcommand{\iy}{\infty}
\newcommand{\fr}{\frac}
\newcommand{\df}{\dfrac}
\newcommand{\pa}{\partial}
\newcommand{\ti}{\times}

\newcommand{\ri}{{\rm i}}

\newcommand{\be}{\begin{equation}}
\newcommand{\ee}{\end{equation}}

\begin{document}

\title[Inverse Electromagnetic Diffraction]{Inverse Electromagnetic Diffraction
by Biperiodic Dielectric Gratings}

\author{Xue Jiang}
\address{School of Science, Beijing University of Posts and
Telecommunications, Beijing 100876, China.}
\email{jxue@lsec.cc.ac.cn}

\author{Peijun Li}
\address{Department of Mathematics, Purdue University, West Lafayette, Indiana
47907, USA.}
\email{lipeijun@math.purdue.edu}

\thanks{The research of XJ was supported in part by China NSF grant 11401040 and
by the Fundamental Research Funds for the Central Universities 24820152015RC17.
The research of PL was partially supported by NSF DMS-1151308.}



\begin{abstract}
Consider the incidence of a time-harmonic electromagnetic plane wave onto a
biperiodic dielectric grating, where the surface is assumed to be a small
and smooth perturbation of a plane. The diffraction is modeled as a transmission
problem for Maxwell's equations in three dimensions. This paper concerns the
inverse diffraction problem which is to reconstruct the grating surface from
either the diffracted field or the transmitted field. A novel approach is
developed to solve the challenging nonlinear and ill-posed inverse problem. The
method requires only a single incident field and is realized via the fast
Fourier transform. Numerical results show that it is simple, fast, and stable to
reconstruct biperiodic dielectric grating surfaces with
super-resolved resolution. 
\end{abstract}

\maketitle

\section{Introduction}

Consider the diffraction of a time-harmonic electromagnetic plane incident wave
by a biperiodic structure, which is called a crossed or two-dimensional grating
in optical community. Given the structure and the incident field, the direct
problem is to determine the diffracted field. The inverse problem is to
reconstruct the grating surface from measured field. This paper concerns the
latter. Diffractive gratings have been widely used in micro-optics including the
design and fabrication of optical elements such as corrective lenses,
anti-reflective interfaces, beam splitters, and sensors. Driven by the
industrial applications, the diffraction grating problems have received
ever-lasting attention in the engineering and applied mathematical communities
\cite{B-D-C, N-S}. An introduction to this topic can be found in the
monograph \cite{Pet}. We refer to \cite{B-C-M} for a comprehensive review on the
mathematical modeling and computational methods for these problems.

The inverse diffraction problems have been studied extensively for
one-dimensional gratings, where the structures are invariant in one direction
and the model of Maxwell's equations can be simplified into the Helmholtz
equation. Mathematical results, such as uniqueness and stability, are
established by many researchers \cite{Bao-1, B-F, B-C-Y, H-K, Kir}.
Computationally, a number of methods are developed \cite{A-K, B-L-W-2, B-L-L,
B-E, E-H-R, Het, I-R}. Numerical solutions can be found in \cite{A-K-Y,
C-G-H-I-R, D-W, K-T} for solving general inverse surface scattering problems.
There are also many work done for two-dimensional gratings, where Maxwell's
equations must be considered. We refer to \cite{Bao-2, B-L-W-1, Dob-3, L-N-1}
for the existence, uniqueness, and numerical approximations of solutions for the
direct problems. Mathematical studies on uniqueness can be found in \cite{Amm,
B-Z-Z, B-Z, H-Y-Z, H-Z, Y-Z} for the inverse problems. Numerical results are
very rare for the inverse problems due to the nonlinearity, ill-posedness,
and large scale computation \cite{L-N-2}. Despite a great number of work done
for the inverse diffraction problems, they addressed the classical inverse
scattering problems. The reconstructed resolution was limited by Rayleigh's
criterion, approximately half of the incident wavelength, also known as the
diffraction limit \cite{Cou}.

Recently, a novel approach has been developed to solve inverse surface
scattering problems in various near-field imaging modalities \cite{B-Li-1,
B-Li-2, B-Li-3, C-L-W} including the inverse electromagnetic diffraction by a
perfectly electrically conducting grating \cite{B-C-L}.  This work presents the
first quantitative method for solving the inverse diffraction problem of
Maxwell's equations with super-resolved resolution. As is known, the perfect
electric conductor is an idealized material exhibiting infinite electrical
conductivity and may not exist in nature. In this paper, we consider a realistic
dielectric grating and the result is closer to practical applications. The more
elaborate techniques differ from the existing work because a complicated
transmission problem of Maxwell's equations needs to be studied. Related work
on near-field imaging may be found in \cite{B-Lin, C-S-1, C-S-2, L-N-2}.

Specifically, we consider the incidence of an electromagnetic plane wave on a
dielectric crossed grating, where the surface is assumed to be a small and
smooth deformation of a plane. The diffraction is modeled as a transmission
problem for Maxwell's equations in three dimensions. The method begins with the
transformed field expansion and reduces the boundary value problem into a
successive sequence of two-point boundary value problems. Dropping higher
order terms in the expansion, we linearize the nonlinear inverse problem and 
obtain explicit reconstruction formulas for both the reflection and
transmission configurations. A spectral cut-off regularization is adopted to
suppress the exponential growth of the noise in the evanescent wave components,
which carry high spatial frequency of the surface and contribute to the super
resolution. The method requires only a single illumination with one
polarization, one frequency, and one incident direction, and is realized via the
fast Fourier transform. The numerical results are computed by using synthetic
scattering data provided by an adaptive edge element method with a perfectly
matched absorbing layer \cite{B-L-W-1}. Two numerical examples, one smooth
surface and one non-smooth surface, are presented to demonstrate the
effectiveness of the proposed method. Careful numerical studies are carried for
the influence of all the parameters on the reconstructions. The results show
that the method is simple, fast, and stable to reconstruct dielectric crossed
grating surfaces with subwavelength resolution.

The paper is organized as follows. In Section 2, the model problem is
introduced. Section 3 presents the transformed field expansion to obtain the 
analytic solution of the direct problem. Explicit reconstruction formulas
are derived for the inverse problem in Section 4. Numerical examples are
reported in Section 5. The paper is concluded with some general remarks and
direction for future work in Section 6.

\section{Model problem}

In this section, we define some notation and introduce a boundary
value problem for the diffraction by a biperiodic dielectric grating.

\subsection{Maxwell's equations}

Let us first specify the problem geometry. Denote $(\rh, z)\in\mb{R}^3$, where
$\rh=(x, y)\in\mb{R}^2$. As seen in Figure \ref{pg}, the problem may be
restricted to a single period of $\La=(\La_1,
\La_2)$ in $\rh$ due to the periodicity of the structure. Let the surface in one
period be described by $S=\{(\rh, z)\in\mb{R}^3: z=\ph(\rh),\,0<x<\La_1,
0<y<\La_2\}$, where $\ph\in C^2(\mb{R}^2)$ is a biperiodic grating surface
function. We assume that
\be\lb{gs}
\ph(\rh)=\de\ps(\rh),
\ee
where $\de>0$ is a small surface deformation parameter, $\ps\in C^2(\mb{R}^2)$
is also a biperiodic function and describes the shape of the
grating surface.

We let $S$ be embedded in the rectangular slab:
\[
 \Om=\{(\rh, z)\in\mb{R}^3: z_-<z<z_+\}=\mb{R}^2\ti(z_-,\,z_+),
\]
where $z_+>0$ and $z_-<0$ are two constants. Hence the domain $\Om$ is bounded
by two plane surfaces $\Ga_\pm=\{(\rh, z)\in\mb{R}^3: z=z_\pm\}$. Let
$\Om_S^+=\{(\rh, z): z>\ph(\rh)\}$ and $\Om_S^-=\{(\rh, z): z<\ph(\rh)\}$ be
filled with homogeneous materials which are characterized by the electric
permittivity $\ve_+$ and $\ve_-$, respectively.

\begin{figure}
\centering
\includegraphics[width=0.4\textwidth]{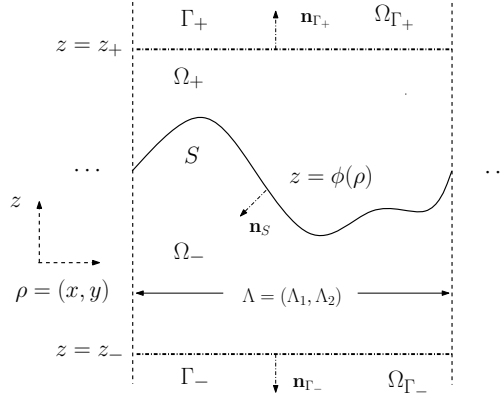}
\caption{The problem geometry of a biperiodic dielectric grating.}
\lb{pg}
\end{figure}

Let $({\bf E}^{\rm inc},\, {\bf H}^{\rm inc})$ be the incoming
electromagnetc plane waves, where
\be\lb{pw}
{\bf E}^{\rm inc}={\bf p} e^{\ri\ka_+(\al\cdot\rh-\ba z)},\qu {\bf H}^{\rm
inc}=\left(\fr{\ve_+}{\mu}\right)^{1/2}{\bf q} e^{\ri\ka_+(\al\cdot\rh-\ba z)}.
\ee
Here $\ka_+=\om(\mu\ve_+)^{1/2}$ is the wavenumber in $\Om_S^+$, $\om>0$ is the
angular frequency, $\mu$ is the magnetic permeability and is assumed to be a
positive constant everywhere, $\al=(\al_1, \al_2)$, $\al_1=\sin\te_1\cos\te_2$,
$\al_2=\sin\te_1\sin\te_2$, and $\ba=\cos\te_1$, where $\te_1$ and $\te_2$ are
the latitudinal and longitudinal incident angles, respectively, which satisfy
$0\le\te_1<\pi/2, 0\le\te_2<2\pi$. Denote by ${\bf d}=(\al_1, \al_2, -\ba)$ the
unit propagation direction vector. The unit polarization vectors ${\bf p}=(p_1,
p_2, p_3)$ and ${\bf q}=(q_1, q_2, q_3)$ satisfy
\[
{\bf p}\cdot{\bf d}=0, \qu {\bf q}={\bf d}\ti{\bf p},
\]
which gives explicitly that 
\[
q_1=\al_2 p_3 + \ba p_2,\qu q_2=-(\al_1 p_3 + \ba p_1),\qu q_3=\al_1 p_2 - \al_2
p_1.
\]
For normal incident, i.e., $\te_1=0$, we have
\[
\al_1=0, \qu \al_2=0, \qu \ba=1,\qu
q_1=p_2, \qu q_2=-p_1, \qu q_3=0.
\]
Hence we get from $|{\bf p}|=|{\bf q}|=1$ that
\[
p_1^2 + p_2^2=1, \qu p_3=0.
\]
For simplicity, we focus on the normal incidence from now on since our method
requires only a single incidence. In fact, this is the most convenient way
to illuminate the grating structure. The method also works for non-normal
incidence with obvious modifications.

Let ${\bf E}^{\rm inc}=(E_1^{\rm inc}, E_2^{\rm inc}, E_3^{\rm inc})$ and
${\bf H}^{\rm inc}=(H_1^{\rm inc}, H_2^{\rm inc}, H_3^{\rm inc})$. Under the
normal incidence, the incoming plane waves \eqref{pw} reduce to
\be\lb{iw}
E_j^{\rm inc}=p_j e^{-\ri\ka_+ z}, \qu H_j^{\rm
inc}=\left(\fr{\ve_+}{\mu}\right)^{1/2} q_j e^{-\ri\ka_+ z},
\ee
which satisfy the time-harmonic Maxwell equation:
\[
\na\ti{\bf E}^{\rm inc}-\ri\om\mu{\bf H}^{\rm inc}=0,\qu \na\ti{\bf H}^{\rm
inc}+\ri\om\ve_+{\bf E}^{\rm inc}=0\qu\text{in}~\Om_S^+.
\]

The time-harmonic electromagnetic waves satisfy Maxwell's equations:
\be\lb{tf}
\na\ti{\bf E}-\ri\om\mu{\bf H}=0, \qu \na\ti{\bf H}+\ri\om\ve{\bf
E}=0\qu\text{in} ~ \mb{R}^3,
\ee
where $({\bf E}, {\bf H})$ are the total electric and magnetic fields, and the
dielectric permittivity
\begin{align*}
\ve=\left\{
\begin{array}{lll}
\ve_+ &\qu\text{in} & \Om_S^+,\\[5pt]
\ve_- &\qu\text{in} & \Om_S^-.
\end{array}
\right.
\end{align*}
Motivated by uniqueness, we are interested in periodic solutions of
$({\bf E}, {\bf H})$ in $\rh$ with period $\La$, i.e., $({\bf E}, {\bf
H})$ satisfy
\[
{\bf E}(\rh+\La, z)={\bf E}(\rh, z),\qu {\bf
H}(\rh+\La, z)={\bf H}(\rh, z).
\]
The total fields can be decomposed into
\begin{align*}
({\bf E},\,{\bf H})=\left\{
\begin{array}{cll}
({\bf E}^{\rm inc},\,{\bf H}^{\rm inc})+({\bf E}^{\rm d},\,{\bf H}^{\rm d})
&\qu\text{in} & \Om_S^+,\\[5pt]
({\bf E}^{\rm t},\,{\bf H}^{\rm t}) &\qu\text{in} & \Om_S^-,
\end{array}
\right.
\end{align*}
where $({\bf E}^{\rm d},\,{\bf H}^{\rm d})$ are the diffracted fields and $({\bf
E}^{\rm t},\,{\bf H}^{\rm t})$ are the transmitted fields. They are required to
satisfy the bounded outgoing wave condition.

\subsection{Transparent boundary condition}

In this section, we introduce transparent boundary conditions on $\Ga_\pm$
which are equivalent to the bounded outgoing wave condition. The detailed
derivation can be found in \cite{B-L-W-2}.

Let $n=(n_1,\,n_2)\in\mb{Z}^2$ and denote $\al_n=(\al_{1 n},\,\al_{2 n})$, where
$\al_{1 n}=2\pi n_1/\La_1$ and $\al_{2 n}=2\pi n_2/\La_2$. For any vector field
${\bf u}=(u_1,\, u_2,\, u_3)$, denote its tangential components on $\Ga_\pm$ by
\[
{\bf u}_{\Ga_\pm}={\bf n}_{\Ga_\pm}\ti({\bf u}\ti{\bf
n}_{\Ga_\pm})=(u_1(\rh, z_\pm),\, u_2(\rh, z_\pm),\,0),
\]
and its tangential traces on $\Ga_\pm$ by
\begin{align*}
 {\bf u}\ti{\bf n}_{\Ga_+}=(u_2(\rh, z_+),\, -u_1(\rh, z_+),\, 0),\\[5pt]
 {\bf u}\ti{\bf n}_{\Ga_-}=(-u_2(\rh, z_-),\, u_1(\rh, z_-),\, 0),
\end{align*}
where ${\bf n}_{\Ga_\pm}=(0,\, 0,\, \pm 1)$ are the unit normal vectors on
$\Ga_\pm$.

For any tangential vector ${\bf u}(\rh, z_+)=(u_1(\rh, z_+),\, u_2(\rh,
z_+),\, 0)$ on $\Ga_+$, where $u_j$ are biperiodic functions in $\rh$ with
period $\La$, we define the capacity operator $T_+$:
\be\lb{bop}
T_+{\bf u}=(v_1(\rh, z_+),\, v_2(\rh, z_+),\, 0),
\ee
where $v_j$ are also biperiodic functions in $\rh$ with the same period
$\La$. Here $u_j$ and $v_j$ have the following Fourier series expansions
\[
u_j(\rh, z_+)=\sum_{n\in\mb{Z}^2} u_{j n}(z_+)e^{\ri\al_n\cdot\rh}, \qu v_j(\rh,
z_+)=\sum_{n\in\mb{Z}^2} v_{j n}(z_+)e^{\ri\al_n\cdot\rh},
\]
and the Fourier coefficients $u_{j n}$ and $v_{j n}$ satisfy
\begin{align*}
\left\{
\begin{array}{c}
v_{1 n}(z_+)=\df{1}{\om\mu\ba^+_n}\left[(\ka_+^2-\al_{2 n}^2)u_{1 n}(z_+)+\al_{1
n}\al_{2 n}u_{2 n}(z_+)\right],\\[12pt]
v_{2 n}(z_+)=\df{1}{\om\mu\ba^+_n}\left[(\ka_+^2-\al_{1 n}^2)u_{2 n}(z_+)+\al_{1
n}\al_{2 n}u_{1 n}(z_+)\right],
\end{array}
\right.
\end{align*}
where
\be\lb{bap}
 (\ba_n^+)^2=\ka^2_+-|\al_n|^2\qu\text{with} ~ {\rm Im}\,\ba_n^+ >0.
\ee
We exclude possible resonance by assuming that $\ba^+_n\ne 0$ for all
$n\in\mb{Z}^2$.

Using the capacity operator \eqref{bop}, we impose a transparent boundary
condition on $\Ga_+$:
\[
T_+({\bf E}_{\Ga_+}-{\bf E}^{\rm inc}_{\Ga_+})=({\bf H}-{\bf H}^{\rm
inc})\ti{\bf n}_{\Ga_+},
\]
which maps the tangential component of the scattered electric field to the
tangential trace of the scattered magnetic field. Equivalently, the above
boundary condition can be written as
\be\lb{tbcp}
(\na\ti{\bf E})\ti{\bf n}_{\Ga_+}=\ri\om\mu T_+{\bf E}_{\Ga_+}+{\bf f},
\ee
where
\[
{\bf f}=\ri\om\mu({\bf H}^{\rm inc}\ti{\bf n}_{\Ga_+}-T_+{\bf E}^{\rm
inc}_{\Ga_+})=(f_1,\, f_2,\, f_3).
\]
Recalling the incident fields \eqref{iw} and using the boundary operator
\eqref{bop}, we have explicitly that
\[
f_j=-2\ri\ka_+ p_j e^{-\ri\ka_+ z_+}.
\]

Similarly, for any given tangential vector ${\bf u}(\rh, z_-)=(u_1(\rh,
z_-),\, u_2(\rh, z_-),\,0)$ on $\Ga_-$, where $u_j(\rh, z_-)$ is a biperiodic
function in $\rh$ with period $\La$, we define the capacity operator $T_-$:
\be\lb{bom}
T_-{\bf u}=(v_1(\rh, z_-),\, v_2(\rh, z_-),\, 0),
\ee
where $v_j$ is also a biperiodic function in $\rh$ with the same period
$\La$. Here $u_j$ and $v_j$ have the following Fourier series expansions
\[
u_j(\rh, z_-)=\sum_{n\in\mb{Z}^2} u_{j n}(z_-)e^{\ri\al_n\cdot\rh}, \qu v_j(\rh,
z_-)=\sum_{n\in\mb{Z}^2} v_{j n}(z_-)e^{\ri\al_n\cdot\rh},
\]
and the Fourier coefficients $u_{j n}$ and $v_{j n}$ satisfy
\begin{align*}
\left\{
\begin{array}{c}
v_{1 n}(z_-)=\df{1}{\om\mu\ba^-_n}\left[(\ka_-^2-\al_{2 n}^2)u_{1 n}(z_-)+\al_{1
n}\al_{2 n}u_{2 n}(z_-)\right],\\[12pt]
v_{2 n}(z_-)=\df{1}{\om\mu\ba^-_n}\left[(\ka_-^2-\al_{1 n}^2)u_{2 n}(z_-)+\al_{1
n}\al_{2 n}u_{1 n}(z_-)\right],
\end{array}
\right.
\end{align*}
where $\ka_-=\om(\mu\ve_-)^{1/2}$ is the wavenumber in $\Om_S^-$ and
\be\lb{bam}
 (\ba_n^-)^2=\ka^2_- -|\al_n|^2\qu\text{with} ~ {\rm Im}\,\ba_n^- >0.
\ee
Here we also assume that $\ba^-_n\ne 0$ for all $n\in\mb{Z}^2$.

Based on \eqref{bom}, a transparent boundary condition may be proposed on
$\Ga_-$:
\[
T_- {\bf E}_{\Ga_-}={\bf H}\ti{\bf n}_{\Ga_-},
\]
which is equivalent to
\be\lb{tbcm}
(\na\ti{\bf E})\ti{\bf n}_{\Ga_-}=\ri\om\mu T_-{\bf E}_{\Ga_-}.
\ee

\subsection{Transmission problem}

Taking curl on both sides of \eqref{tf}, we may eliminate the magnetic field and
obtain a decoupled equation for the electric field:
\be\lb{me}
\na\ti(\na\ti{\bf E})-\ka^2{\bf E}=0\qu\text{in} ~ \Om,
\ee
where the wavenumber
\begin{align*}
\ka=\left\{
\begin{array}{lll}
\ka_+ &\qu\text{in} & \Om_S^+,\\[5pt]
\ka_- &\qu\text{in} & \Om_S^-.
\end{array}
\right.
\end{align*}

Denote $\Om_+=\Om_S^+\cap\Om=\{(\rh, z): \ph(\rh)<z<z_+\}$ and
$\Om_-=\Om_S^-\cap\Om=\{(\rh, z): z_-<z<\ph(\rh)\}$. Let ${\bf E}^+$ and ${\bf
E}^-$ be the restriction of ${\bf E}$ in $\Om_+$ and $\Om_-$, respectively,
i.e., ${\bf E}^\pm={\bf E}|_{\Om_\pm}$. It is useful to have an
equivalent scalar form of \eqref{me} when applying the transformed field
expansion. Denote ${\bf E}^\pm=(E_1^\pm,\, E_2^\pm,\, E_3^\pm)$. We may
reformulate \eqref{me} into the Helmholtz equation:
\be\lb{ctf}
\De E_j^\pm+\ka^2_\pm E_j^\pm=0\qu\text{in} ~\Om_\pm.
\ee

The transparent boundary conditions \eqref{tbcp} and \eqref{tbcm} can be written
as
\begin{align}\lb{ctbcp}
\left\{
\begin{array}{c}
\pa_z E_1^+ -\pa_x E_3^+=\ri\om\mu H_1^+ +f_1,\\[5pt]
\pa_z E_2^+ -\pa_y E_3^+=\ri\om\mu H_2^+ +f_2,
\end{array}
\right.
\end{align}
and
\begin{align}\lb{ctbcm}
\left\{
\begin{array}{c}
\pa_z E_1^- -\pa_x E_3^-=-\ri\om\mu H_1^-,\\[5pt]
\pa_z E_2^- -\pa_y E_3^-=-\ri\om\mu H_2^-,
\end{array}
\right.
\end{align}
where the Fourier coefficients of the periodic functions $H_1^\pm$ and $H_2^\pm$
are given by
\begin{align*}
\left\{
\begin{array}{c}
H_{1 n}^\pm(z_\pm)=\df{1}{\om\mu\ba^\pm_n}\left[(\ka^2_\pm-\al_{2 n}^2) E_{1
n}^\pm(z_\pm)+\al_{1 n}\al_{2 n}E_{2 n}^\pm(z_\pm)\right],\\[12pt]
H_{2 n}^\pm(z_\pm)=\df{1}{\om\mu\ba^\pm_n}\left[(\ka^2_\pm-\al_{1 n}^2) E_{2
n}^\pm(z_\pm)+\al_{1 n}\al_{2 n}E_{1 n}^\pm(z_\pm)\right].
\end{array}
\right.
\end{align*}
Here $E_{1 n}^\pm(z_\pm)$ and $E_{2 n}^\pm(z_\pm)$ are the Fourier coefficients
of the periodic electric field $E_1^\pm(\rh, z_\pm)$ and $E_2^\pm(\rh,
z_\pm)$, respectively.

The continuity conditions are needed to reformulate the boundary value problem
into a transmission problem. It is known that the tangential traces of
the electric and magnetic fields are continuous across the grating surface,
i.e.,
\[
 {\bf E}^+\ti{\bf n}_S ={\bf E}^-\ti{\bf n}_S,\qu {\bf H}^+\ti{\bf n}_S={\bf
H}^-\ti{\bf n}_S,\qu z=\ph(\rh),
\]
where ${\bf n}_S=(\ph_x, \ph_y, -1)$ is the normal vector on $S$ pointing
from $\Om_S^+$ to $\Om_S^-$. Explicitly, we have the continuity conditions
\begin{align}\lb{cce}
\left\{
\begin{array}{c}
E_2^+ +\ph_y E_3^+ =E_2^- +\ph_y E_3^-,\\[5pt]
E_1^+ +\ph_x E_3^+ =E_1^- +\ph_x E_3^-,
\end{array}
\right.
\end{align}
and
\begin{align}\lb{cch}
\left\{
\begin{array}{l}
 \left(\pa_z E_1^+ -\pa_x E_3^+\right)+\ph_y\left(\pa_x E_2^+
-\pa_y E_1^+\right)\\[5pt]
\hs{2cm}=(\pa_z E_1^- -\pa_x E_3^-)+\ph_y(\pa_x E_2^-
-\pa_y E_1^-),\\[10pt]
\left(\pa_y E_3^+ -\pa_z E_2^+\right)+\ph_x\left(\pa_x E_2^+
-\pa_y E_1^+\right)\\[5pt]
\hs{2cm}=\left(\pa_y E_3^- -\pa_z E_2^-\right)+\ph_x\left(\pa_x
E_2^- -\pa_y E_1^-\right).
\end{array}
\right.
\end{align}

The transparent boundary conditions \eqref{ctbcp}, \eqref{ctbcm} and the
continuity conditions \eqref{cce} and \eqref{cch} are not enough to determine
the fields ${\bf E}^\pm_j$. Additional information can be obtained from the
divergence free conditions
\be\lb{cdf}
\pa_x E_1^\pm +\pa_y E_2^\pm +\pa_z E_3^\pm=0\qu\text{in} ~ \Om_\pm.
\ee

Given the grating surface function $\ph(\rh)$, the direct problem is to
determine the fields ${\bf E}^\pm_j$. This work is focused on the inverse
problem, which is to reconstruct the grating surface function $\ph(\rh)$ from
the tangential traces of the total field measured at either $\Ga_+$, i.e., ${\bf
E}(\rh, z_+)\ti{\bf n}_{\Ga_+}=(E_1(\rh, z_+),\, E_2(\rh, z_+),\,0)$ called the
reflection configuration, or $\Ga_-$, i.e., ${\bf E}(\rh, z_-)\ti{\bf
n}_{\Ga_-}=(-E_1(\rh, z_-),\, E_2(\rh, z_-),\,0)$ called the transmission
configuration. In particular, we are interested in the inverse problem in
near-field regime where the measurement distance $|z_\pm|$ is much smaller than
the wavelength $\la=2\pi/\kappa_+$ of the incident field.

\section{Transformed field expansion}

In this section, we introduce the transformed field expansion to analytically
derive the solution for the direct problem. We refer to \cite{B-R-1, N-R-1} for
solving the direct surface scattering problems by using the transformed field
expansion and related boundary perturbation method.

\subsection{Change of variables}

Consider the change of variables: 
\[
\tl{x}=x, \qu \tl{y}=y, \qu \tl{z}=z_+\left(\fr{z-\ph}{z_+ -\ph}\right),\qu
\ph<z<z_+,
\]
and
\[
\tl{x}=x, \qu \tl{y}=y, \qu \tl{z}=z_-\left(\fr{z-\ph}{z_- -\ph}\right),\qu
z_-<z<\ph,
\]
which maps the domain $\Om_+$ and $\Om_-$ into rectangular slabs
$D_+=\{(\tl{\rh}, \tl{z})\in\mb{R}^3: 0<\tl{z}<z_+\}$ and $D_-=\{(\tl{\rh},
\tl{z})\in\mb{R}^3: z_-<z<0\}$, respectively.

We seek to restate the diffractive grating problem in the new coordinate.
Introduce a new function $\tl{\bf E}^\pm=(\tl{E}_1^\pm,\, \tl{E}_2^\pm,\,
\tl{E}_3^\pm)$ and let $\tl{E}_j^\pm(\tl{x},\, \tl{y},\, \tl{z})=E_j^\pm(x,\,
y,\, z)$ under the transformation. After tedious but straightforward
calculations, it can be verified from \eqref{ctf} that the total electric field,
upon dropping the tilde, satisfies the equation
\begin{align}\lb{ttf}
c_1^\pm\fr{\pa^2 E_j^\pm}{\pa x^2}&+c_1^\pm\fr{\pa^2 E_j^\pm}{\pa y^2}+
c_2^\pm \fr{\pa^2 E_j^\pm}{\pa z^2}-c_3^\pm \fr{\pa^2 E_j^\pm}{\pa x\pa
z}\nt\\[5pt]
&-c_4^\pm \fr{\pa^2 E_j^\pm}{\pa y\pa z}-c_5^\pm\fr{\pa E_j^\pm}{\pa
z}+\ka^2_\pm c_1^\pm E_j^\pm=0\qu\text{in}~ D_\pm,
\end{align}
where
\begin{align*}
c_1^\pm&=(z_\pm-\ph)^2,\\
c_2^\pm&=(\ph_x^2+\ph_y^2)(z_\pm-z)^2+z_\pm^2,\\
c_3^\pm&=2\ph_x(z_\pm-z)(z_\pm-\ph),\\
c_4^\pm&=2\ph_y(z_\pm-z)(z_\pm-\ph),\\
c_5^\pm&=(z_\pm-z)\big[(\ph_{xx}+\ph_{yy})(z_\pm-\ph)+2(\ph_x^2+\ph_y^2)\big].
\end{align*}
The transparent boundary conditions \eqref{ctbcp} and \eqref{ctbcm} reduce to
\begin{align}\lb{ttbcp}
\left\{
\begin{array}{c}
\left(\df{z_+}{z_+ -\ph}\right)\pa_z E_1^+ -\pa_x E_3^+=\ri\om\mu
H_1^+ +f_1,\\[12pt]
\left(\df{z_+}{z_+ -\ph}\right)\pa_z E_2^+ -\pa_y E_3^+=\ri\om\mu
H_2^+ +f_2,
\end{array}
\right.
\end{align}
and
\begin{align}\lb{ttbcm}
\left\{
\begin{array}{c}
\left(\df{z_-}{z_- -\ph}\right)\pa_z E_1^- -\pa_x E_3^- =-\ri\om\mu
H_1^-,\\[12pt]
\left(\df{z_-}{z_- -\ph}\right)\pa_z E_2^- -\pa_y E_3^- =-\ri\om\mu
H_2^-.
\end{array}
\right.
\end{align}
The continuity conditions \eqref{cce} and \eqref{cch} are changed to
\begin{align}\lb{tce}
\left\{
\begin{array}{c}
E_2^+ +\ph_y E_3^+ =E_2^- +\ph_y E_3^-,\\[5pt]
E_1^+ +\ph_x E_3^+ =E_1^- +\ph_x E_3^-,
\end{array}
\right.
\end{align}
and
\begin{align}\lb{tch}
\left\{
\begin{array}{r}
 \left(\df{z_+}{z_+-\ph}\right)\Bigl[\ph_x \pa_z E_3^+ -\ph_x\ph_y \pa_z
E_2^+ + (1+\ph_y^2)\pa_z E_1^+\Bigr]\\[5pt]
 -\left(\pa_x E_3^+ -\ph_y \pa_x E_2^+ +\ph_y\pa_y E_1^+\right) \\
 = \left(\df{z_-}{z_--\ph}\right)\Bigl[\ph_x \pa_z E_3^-
-\ph_x\ph_y \pa_z E_2^- + (1+\ph_y^2)\pa_z E_1^-\Bigr]\\[5pt]
-\left(\pa_x E_3^- -\ph_y \pa_x E_2^- +\ph_y\pa_y E_1^-\right),\\[10pt]
\left(\df{z_+}{z_+-\ph}\right)\Bigl[\ph_y \pa_z E_3^+ + (1+\ph_x^2)\pa_z E_2^+
-\ph_x\ph_y\pa_z E_1^+\Bigr] \\[5pt]
-\left(\pa_y E_3^+ +\ph_x\pa_x E_2^+ -\ph_x\pa_y E_1^+\right) \\
= \left(\df{z_-}{z_- -\ph}\right)\Bigl[\ph_y \pa_z E_3^- +
(1+\ph_x^2)\pa_z E_2^- -\ph_x\ph_y\pa_z E_1^-\Bigr]\\[5pt]
 -\left(\pa_y E_3^- +\ph_x\pa_x E_2^- -\ph_x\pa_y E_1^-\right) .
\end{array}
\right.
\end{align}
The divergence free condition \eqref{cdf} becomes
\begin{align}\lb{tdf}
\pa_x E^\pm_1+\pa_y E^\pm_2 -\left(\fr{z_\pm-z}{z_\pm-\ph}\right)(\ph_x\pa_z
E^\pm_1+\ph_y\pa_z E^\pm_2)\nt\\[5pt]
+\left(\fr{z_\pm}{z_\pm-\ph}\right)\pa_z
E^\pm_3=0\qu\text{in}~ D_\pm.
\end{align}

\subsection{Power series}

Recalling $\ph=\de\ps$ in \eqref{gs}, we use a classical boundary perturbation
argument and consider a formal expansion of $E_j^\pm$ in a power series of
$\de$:
\be\lb{ps}
E_j^\pm(\rh, z; \de)=\sum_{k=0}^\iy E^{\pm (k)}_j(\rh, z)\, \de^k.
\ee

Substituting $\ph=\de\ps$ and the power series expansion \eqref{ps} into
$c_j^\pm$ and \eqref{ttf}, we may derive a recursion equation for $E^{\pm
(k)}_j$:
\be\lb{rtf}
\De E_j^{\pm (k)}+\ka^2_\pm E^{\pm (k)}_j=F^{\pm (k)}_j\qu\text{in} ~ D_\pm,
\ee
where the nonhomogeneous term
\begin{align*}
F^{\pm (k)}_j&=\fr{2\ps}{z_\pm}\fr{\pa^2 E_j^{\pm (k-1)}}{\pa
x^2}+\fr{2\ps}{z_\pm}\fr{\pa^2 E_j^{\pm (k-1)}}{\pa
y^2}+\fr{2(z_\pm -z)\ps_x}{z_\pm}\fr{\pa^2 E_j^{\pm (k-1)}}{\pa x\pa
z}\\[5pt]
&+\fr{2(z_\pm-z)\ps_y}{z_\pm}\fr{\pa^2 E_j^{\pm (k-1)}}{\pa y\pa z}
+\fr{(z_\pm -z)(\ps_{xx}+\ps_{yy})}{z_\pm}\fr{\pa E_j^{\pm (k-1)}}{\pa
z}+\fr{2\ka^2_\pm\ps}{z_\pm}E_j^{\pm (k-1)}\\[5pt]
&-\fr{\ps^2}{z_\pm^2}\fr{\pa^2
E_j^{\pm (k-2)}}{\pa
x^2}-\fr{\ps^2}{z_\pm^2}\fr{\pa^2 E_j^{\pm (k-2)}}{\pa
y^2}-\fr{(z_\pm-z)^2(\ps_x^2+\ps_y^2)}{z_\pm^2}\fr{\pa^2 E_j^{\pm (k-2)}}{\pa
z^2}\\[5pt]
&-\fr{2\ps\ps_x(z_\pm -z)}{z_\pm^2}\fr{\pa^2 E_j^{\pm (k-2)}}{\pa x\pa
z}-\fr{2\ps\ps_y(z_\pm -z)}{z_\pm^2}\fr{\pa^2 E_j^{\pm (k-2)}}{\pa y\pa
z}\\[5pt]
&+\fr{(z_\pm
-z)\big[2(\ps_x^2+\ps_y^2)-\ps(\ps_{xx}+\ps_{yy})\big]}{z_\pm^2}\fr{\pa
E_j^{\pm (k-2)}}{\pa z}-\fr{\ka_\pm^2\ps^2}{z_\pm^2}E_j^{\pm (k-2)}.
\end{align*}
Here $\ps_x=\pa_x\ps(x, y)$ and $\ps_y=\pa_y\ps(x, y)$ are the partial
derivatives.

Substituting \eqref{ps} into the transparent boundary conditions \eqref{ttbcp}
and \eqref{ttbcm}, we obtain
\begin{align}\lb{rtbcp}
\left\{
\begin{array}{c}
\pa_z E_1^{+ (k)}-\pa_x E_3^{+ (k)}=\ri\om\mu H_1^{+ (k)}-f_1^{+ (k)},\\[5pt]
\pa_z E_2^{+ (k)}-\pa_y E_3^{+ (k)}=\ri\om\mu H_2^{+ (k)}-f_2^{+ (k)},
\end{array}
\right.
\end{align}
and
\begin{align}\lb{rtbcm}
\left\{
\begin{array}{c}
\pa_z E_1^{- (k)}-\pa_x E_3^{- (k)}=-\ri\om\mu H_1^{- (k)}-f_1^{- (k)},\\[5pt]
\pa_z E_2^{- (k)}-\pa_y E_3^{- (k)}=-\ri\om\mu H_2^{- (k)}-f_2^{- (k)},
\end{array}
\right.
\end{align}
where
\begin{align*}
f_1^{+ (0)}=-f_1, ~ f_1^{+ (1)}=\fr{\ps}{z_+}\pa_z E_1^{+ (0)}, ~ f_1^{+
(k)}=\fr{\ps}{z_+}\left(\pa_x E_3^{+ (k-1)}+\ri\om\mu H_1^{+
(k-1)}\right),\\[5pt]
f_2^{+ (0)}=-f_2, ~ f_2^{+ (1)}=\fr{\ps}{z_+}\pa_z E_2^{+ (1)}, ~ f_2^{+
(k)}=\fr{\ps}{z_+}\left(\pa_y E_3^{+ (k-1)}+\ri\om\mu H_2^{+ (k-1)}\right),
\end{align*}
and
\begin{align*}
f_1^{- (0)}=0, ~ f_1^{- (1)}=\fr{\ps}{z_-}\pa_z E_1^{- (0)}, ~ f_1^{-
(k)}=\fr{\ps}{z_-}\left(\pa_x E_3^{- (k-1)}-\ri\om\mu H_1^{-
(k-1)}\right),\\[5pt]
f_2^{- (0)}=0, ~ f_2^{- (1)}=\fr{\ps}{z_-}\pa_z E_2^{- (0)}, ~ f_2^{-
(k)}=\fr{\ps}{z_-}\left(\pa_y E_3^{- (k-1)}-\ri\om\mu H_2^{- (k-1)}\right).
\end{align*}
Here the Fourier coefficients of $H^{\pm (k)}_1(\rh, z_\pm)$ and $H^{\pm
(k)}_2(\rh, z_\pm)$ are
\begin{align*}
\left\{
\begin{array}{c}
H^{\pm (k)}_{1 n}(z_\pm)=\df{1}{\om\mu\ba^\pm_n}\left[(\ka_\pm^2-\al_{2
n}^2)E^{\pm (k)}_{1 n}(z_\pm)+\al_{1 n}\al_{2 n}E^{\pm
(k)}_{2 n}(z_\pm)\right],\\[12pt]
H^{\pm (k)}_{2 n}(z_\pm)=\df{1}{\om\mu\ba^\pm_n}\left[(\ka_\pm^2-\al_{1
n}^2)E^{\pm (k)}_{2 n}(z_\pm)+\al_{1 n}\al_{2 n}E^{\pm (k)}_{1 n}(z_\pm)\right].
\end{array}
\right.
\end{align*}
Here $E^{\pm (k)}_{1 n}(z_\pm)$ and $E^{\pm (k)}_{2 n}(z_\pm)$ are the Fourier
coefficients of $E^{\pm (k)}_1(\rh, z_\pm)$ and $E^{\pm (k)}_1(\rh,
z_\pm)$, respectively.

Plugging \eqref{ps} into the jump conditions \eqref{tce} and \eqref{tch}
yields
\begin{align}\lb{rce}
\left\{
\begin{array}{c}
E_2^{+ (k)} +\ps_y E_3^{+ (k-1)} =E_2^{- (k)} +\ps_y E_3^{- (k-1)},\\[5pt]
E_1^{+ (k)} +\ps_x E_3^{+ (k-1)} =E_1^{- (k)} +\ps_x E_3^{- (k-1)},
\end{array}
\right.
\end{align}
and
\begin{align}\lb{rch}
\left\{
\begin{array}{r}
 \left(\pa_z E_1^{+ (k)}+\ps_x \pa_z E_3^{+ (k-1)} +\ps_y^2\pa_z
E_1^{+ (k-2)}-\ps_x\ps_y \pa_z E_2^{+ (k-2)}\right)\\[5pt]
-z_-^{-1}\Bigl(\pa_z E_1^{+ (k-1)}+\ps_x \pa_z E_3^{+ (k-2)}
+\ps_y^2\pa_z E_1^{+ (k-3)}-\ps_x\ps_y \pa_z E_2^{+ (k-3)}\Bigr)\ps\\[5pt]
-\left(\pa_x E_3^{+ (k)}+\ps_y \pa_y E_1^{+ (k-1)}-\ps_y \pa_x E_2^{+
(k-1)}\right)\\[5pt]
+\left(z_+^{-1} + z_-^{-1}\right)\left(\pa_x E_3^{+ (k-1)}+\ps_y \pa_y
E_1^{+ (k-2)}-\ps_y\pa_x E_2^{+ (k-2)}\right)\ps\\[5pt]
-(z_+ z_-)^{-1}\left(\pa_x E_3^{+ (k-2)}+\ps_y\pa_y E_1^{+
(k-3)}-\ps_y\pa_x E_2^{+ (k-3)}\right)\ps^2\\[5pt]
=\left(\pa_z E_1^{- (k)}+\ps_x \pa_z E_3^{- (k-1)} +\ps_y^2\pa_z
E_1^{- (k-2)}-\ps_x\ps_y \pa_z E_2^{- (k-2)}\right)\\[5pt]
-z_+^{-1}\left(\pa_z E_1^{- (k-1)}+\ps_x \pa_z E_3^{- (k-2)}
+\ps_y^2\pa_z E_1^{- (k-3)}-\ps_x\ps_y \pa_z E_2^{- (k-3)}\right)\ps\\[5pt]
-\left(\pa_x E_3^{- (k)}+\ps_y \pa_y E_1^{- (k-1)}-\ps_y \pa_x E_2^{-
(k-1)}\right)\\[5pt]
+(z_+^{-1} + z_-^{-1})\left(\pa_x E_3^{- (k-1)}+\ps_y \pa_y E_1^{-
(k-2)}-\ps_y\pa_x E_2^{- (k-2)}\right)\ps \\[5pt]
-(z_+ z_-)^{-1} \left(\pa_x E_3^{- (k-2)}+\ps_y\pa_y E_1^{-
(k-3)}-\ps_y \pa_x E_2^{- (k-3)}\right)\ps^2\\[10pt]
 \left(\pa_z E_2^{+ (k)}+\ps_y \pa_z E_3^{+ (k-1)} +\ps_x^2\pa_z
E_2^{+ (k-2)}-\ps_x\ps_y \pa_z E_1^{+ (k-2)}\right)\\[5pt]
-z_-^{-1}\left(\pa_z E_2^{+ (k-1)}+\ps_y \pa_z E_3^{+ (k-2)}
+\ps_x^2\pa_z E_2^{+ (k-3)}-\ps_x\ps_y \pa_z E_1^{+ (k-3)}\right)\ps\\[5pt]
-\left(\pa_y E_3^{+ (k)}+\ps_x (\pa_x E_2^{+ (k-1)}-\pa_y E_1^{+ (k-1)})\right)
\\[5pt]
+(z_+^{-1} + z_-^{-1})\left(\pa_y E_3^{+ (k-1)}+\ps_x\pa_x E_2^{+
(k-2)}-\ps_x\pa_y E_1^{+ (k-2)}\right)\ps\\[5pt]
-(z_+ z_-)^{-1}\left(\pa_y E_3^{+ (k-2)}+\ps_x\pa_x E_2^{+
(k-3)}-\ps_x\pa_y E_1^{+ (k-3)}\right)\ps^2\\[5pt]
=\left(\pa_z E_2^{- (k)}+\ps_y \pa_z E_3^{- (k-1)} +\ps_x^2\pa_z
E_2^{- (k-2)}-\ps_x\ps_y \pa_z E_1^{- (k-2)}\right)\\[5pt]
-z_+^{-1}\left(\pa_z E_2^{- (k-1)}+\ps_y \pa_z E_3^{- (k-2)}
+\ps_x^2\pa_z E_2^{- (k-3)}-\ps_x\ps_y \pa_z E_1^{- (k-3)}\right)\ps\\[5pt]
-\left(\pa_y E_3^{- (k)}+\ps_x\pa_x E_2^{- (k-1)}-\ps_x\pa_y E_1^{-
(k-1)}\right)\\[5pt]
+(z_+^{-1} + z_-^{-1})\left(\pa_y E_3^{- (k-1)}+\ps_x\pa_x E_2^{-
(k-2)}-\ps_x\pa_y E_1^{- (k-2)}\right)\ps \\[5pt]
-(z_+ z_-)^{-1}\left(\pa_y E_3^{- (k-2)}+\ps_x\pa_x E_2^{-
(k-3)}-\ps_x\pa_y E_1^{- (k-3)}\right)\ps^2.
\end{array}
\right.
\end{align}

Substituting \eqref{ps} into the divergence free condition \eqref{tdf} yields
\be\lb{rdf}
\pa_x E_1^{\pm (k)}+\pa_y E_2^{\pm (k)}+\pa_z E_3^{\pm
(k)}=g^{\pm (k)}\qu\text{in}~ D_\pm,
\ee
where
\begin{align*}
w^{\pm (k)}=&\fr{\ps}{z_\pm}\left(\pa_x E_1^{\pm (k-1)}+\pa_y
E_2^{\pm (k-1)}\right)\\[5pt]
&+\left(\fr{z_\pm-z}{z_\pm}\right)\left(\ps_x \pa_z E_1^{\pm (k-1)}+\ps_y
\pa_z E_2^{\pm (k-1)}\right).
\end{align*}

\subsection{Zeroth order}

Recalling the recurrence relation \eqref{rtf} and letting $k=0$, we have
\be\lb{rtf_0}
\De E_j^{\pm (0)}+\ka_\pm^2 E_j^{\pm (0)}=0\qu\text{in} ~ D_\pm.
\ee
The transparent boundary conditions \eqref{rtbcp} and \eqref{rtbcm} becomes
\begin{align}\lb{rtbcp_0}
\left\{
\begin{array}{c}
\pa_z E_1^{+ (0)}(\rh, z_+)-\pa_x E_3^{+ (0)}(\rh, z_+)=\ri\om\mu
H_1^{+ (0)}(\rh, z_+) + f_1(\rh),\\[5pt]
\pa_z E_2^{+ (0)}(\rh, z_+)-\pa_y E_3^{+ (0)}(\rh, z_+)=\ri\om\mu
H_2^{+ (0)}(\rh, z_+)+f_2(\rh),
\end{array}
\right.
\end{align}
and
\begin{align}\lb{rtbcm_0}
\left\{
\begin{array}{c}
\pa_z E_1^{- (0)}(\rh, z_-)-\pa_x E_3^{- (0)}(\rh, z_-)=-\ri\om\mu
H_1^{- (0)}(\rh, z_-),\\[5pt]
\pa_z E_2^{- (0)}(\rh, z_-)-\pa_y E_3^{- (0)}(\rh, z_-)=-\ri\om\mu
H_2^{- (0)}(\rh, z_-).
\end{array}
\right.
\end{align}
The jump conditions \eqref{rce} and \eqref{rch} reduce to
\begin{equation}\lb{rce_0}
E_2^{+ (0)}(\rh, 0) =E_2^{- (0)}(\rh, 0),\qu E_1^{+ (0)}(\rh, 0) =E_1^{-
(0)}(\rh, 0),
\end{equation}
and
\begin{align}\lb{rch_0}
\left\{
\begin{array}{c}
\pa_z E_1^{+ (0)}(\rh, 0) - \pa_x E_3^{+ (0)}(\rh, 0)=\pa_z E_1^{- (0)}(\rh, 0)
- \pa_x E_3^{- (0)}(\rh, 0),\\[5pt]
\pa_z E_2^{+ (0)}(\rh, 0) - \pa_y E_3^{+ (0)}(\rh, 0)=\pa_z E_2^{- (0)}(\rh, 0)
- \pa_y E_3^{- (0)}(\rh, 0).
\end{array}
\right.
\end{align}
The divergence free condition \eqref{rdf} reduces to
\be\lb{rdf_0}
\pa_x E_1^{\pm (0)}+\pa_y E_2^{\pm (0)}+\pa_z E_3^{\pm (0)}=0\qu\text{in} ~
D_\pm.
\ee

Since $E^{\pm (0)}_j(\rh, z)$ and $f_j$ are periodic functions of $\rh$ with
period $\La$, they have the following Fourier expansion
\be\lb{E_0f}
E^{\pm (0)}_j(\rh, z)=\sum_{n\in\mb{Z}^2}E^{\pm (0)}_{j n}(z)
e^{\ri\al_n\cdot\rh}, \qu f_j(\rh)=\sum_{n\in\mb{Z}^2} f_{j
n}e^{\ri\al_n\cdot\rh},
\ee
where
\begin{align*}
f_{j n}=\left\{
\begin{array}{cll}
-2\ri\ka_+ p_j e^{-\ri\ka_+ z_+} &\qu\text{for} & n=0,\\[5pt]
0 &\qu\text{for} & n\ne 0.
\end{array}
\right.
\end{align*}
Plugging \eqref{E_0f} into \eqref{rtf_0}--\eqref{rdf_0}, we obtain an ordinary
differential equation
\be\lb{ode_0f}
\fr{{\rm d}^2 E^{\pm (0)}_{j n}(z)}{{\rm
d}z^2}+(\ba^\pm_n)^2 E^{\pm (0)}_{j n}(z)=0,
\ee
together with the boundary conditions at $z=z_+$:
\begin{align}\lb{rbcp_0f}
\left\{
\begin{array}{l}
{E_{1 n}^{+ (0)}}'-\ri\al_{1 n} E_{3
n}^{+ (0)}=\df{\ri}{\ba^+_n}\left[(\ka_+^2-\al_{2
n}^2)E^{+ (0)}_{1 n}+\al_{1 n}\al_{2 n}E^{+ (0)}_{2 n}\right]+f_{1
n},\\[12pt]
{E_{2 n}^{+ (0)}}'-\ri\al_{2 n} E_{3
n}^{+ (0)}=\df{\ri}{\ba^+_n}\left[(\ka_+^2-\al_{1
n}^2)E^{+ (0)}_{2 n}+\al_{1 n}\al_{2 n}E^{+ (0)}_{1 n}\right] +f_{2
n},\\[12pt]
{E_{3 n}^{+ (0)}}'+\ri\al_{1 n}E_{1 n}^{+ (0)}+\ri\al_{2 n}E_{2
n}^{+ (0)}=0.
\end{array}
\right.
\end{align}
and the boundary conditions at $z=z_-$:
\begin{align}\lb{rbcm_0f}
\left\{
\begin{array}{l}
{E_{1 n}^{- (0)}}'-\ri\al_{1 n} E_{3
n}^{- (0)}=-\df{\ri}{\ba^-_n}\left[(\ka_-^2-\al_{2
n}^2)E^{- (0)}_{1 n}+\al_{1 n}\al_{2 n}E^{- (0)}_{2 n}\right],\\[12pt]
{E_{2 n}^{- (0)}}'-\ri\al_{2 n} E_{3
n}^{- (0)}=-\df{\ri}{\ba^-_n}\left[(\ka_-^2-\al_{1
n}^2)E^{- (0)}_{2 n}+\al_{1 n}\al_{2 n}E^{- (0)}_{1 n}\right],\\[12pt]
{E_{3 n}^{- (0)}}'+\ri\al_{1 n}E_{1 n}^{- (0)}+\ri\al_{2 n}E_{2
n}^{- (0)}=0,
\end{array}
\right.
\end{align}
and the jump conditions at $z=0$:
\begin{equation}\lb{rce_0f}
E_{2 n}^{+ (0)} =E_{2 n}^{- (0)},\qu E_{1 n}^{+ (0)} =E_{1 n}^{- (0)},
\end{equation}
and
\begin{align}\lb{rch_0f}
\left\{
\begin{array}{c}
 {E_{1 n}^{+ (0)}}' -\ri\al_{1 n} E_{3 n}^{+
(0)}={E_{1 n}^{- (0)}}' -\ri\al_{1 n} E_{3 n}^{-
(0)},\\[5pt]
{E_{2 n}^{+ (0)}}' - \ri\al_{2 n} E_{3 n}^{+ (0)}={E_{2 n}^{-
(0)}}' - \ri\al_{2 n}E_{3 n}^{- (0)}.
\end{array}
\right.
\end{align}

It can be verified that the general solutions of the homogeneous second order
equations \eqref{ode_0f} are
\be\lb{sode_0f}
 E_{j n}^{\pm (0)}(z)= A^\pm_{j n} e^{\ri\ba_n^\pm z} + B^\pm_{j
n}e^{-\ri\ba_n^\pm z},
\ee
where $A^\pm_{j n}, B^\pm_{j n}\in\mb{C}$ are to be determined. Substituting
\eqref{sode_0f} into the boundary conditions \eqref{rbcp_0f} and
\eqref{rbcm_0f}, and the jump conditions \eqref{rce_0f} and \eqref{rch_0f}, we
may deduce that
\[
A^{+ (0)}_{j 0}=r p_j,\qu B^{+ (0)}_{j 0}=p_j,\qu A^{- (0)}_{j 0}=0,\qu B^{-
(0)}_{j 0}=t p_j,
\]
and $A^\pm_{j n}=B^\pm_{j n}=A^\pm_{j n}=B^\pm_{j n}=0$ for $n\ne 0$, where
\[
r=\fr{\ka_+ -\ka_-}{\ka_+ + \ka_-}\qu\text{and}\qu
t=\fr{2\ka_+}{\ka_+
+ \ka_-}
\]
are known as the reflection coefficient and the transmission coefficient,
respectively. Hence we find the analytic expressions for the zeroth
order terms:
\begin{align}\lb{s0}
\left\{
\begin{array}{l}
E_j^{+ (0)}(\rh, z)=p_j ( e^{-\ri\ka_+ z}+r e^{\ri\ka_+ z}),\\[10pt]
E_j^{- (0)}(\rh, z)=p_j t e^{-\ri\ka_- z}.
\end{array}
\right.
\end{align}
Clearly, the zeroth order terms consist of the incident wave, the reflected
wave, and the transmitted wave, which come from the diffraction of an
electromagnetic plane wave by a planar surface.

\subsection{First order}

In this section, we derive analytic expressions of the first order terms, and
particularly a connection between their Fourier coefficients and the Fourier
coefficient of the grating profile.

Taking $k=1$ in \eqref{rtf} yields
\be\lb{rtf_1}
\De E_j^{\pm (1)}+\ka^2 E_j^{\pm (1)}=F_j^{\pm (1)}\qu\text{in} ~ D_\pm,
\ee
where
\begin{align*}
F^{\pm (1)}_j&=\fr{2\ps}{z_\pm}\fr{\pa^2 E_j^{\pm (0)}}{\pa
x^2}+\fr{2\ps}{z_\pm}\fr{\pa^2 E_j^{\pm (0)}}{\pa
y^2}+\fr{2(z_\pm -z)\ps_x}{z_\pm}\fr{\pa^2 E_j^{\pm (0)}}{\pa x\pa
z}+\fr{2(z_\pm -z)\ps_y}{z_\pm}\fr{\pa^2 E_j^{\pm (0)}}{\pa y\pa z}\\
&+\fr{(z_\pm -z)(\ps_{xx}+\ps_{yy})}{z_\pm}\fr{\pa E_j^{\pm (0)}}{\pa
z}+\fr{2\ka_\pm^2\ps}{z_\pm}E_j^{\pm (0)}.
\end{align*}
It follows from the explicit expression of the zeroth order term \eqref{s0} that
we have
\begin{align*}
 F_j^{+ (1)}(\rh, z)&=\fr{2\ka^2_+ p_j}{z_+}\left(e^{-\ri\ka_+
z}+ re^{\ri\ka_+ z}\right)\ps\\
&-\fr{\ri\ka_+ p_j(z_+ -z)}{z_+}\left(e^{-\ri\ka_+
z}-re^{\ri\ka_+ z}\right)(\ps_{xx}+\ps_{yy})
\end{align*}
and
\[
 F_j^{- (1)}(\rh, z)=\fr{2\ka_-^2 p_j}{z_-}t e^{-\ri\ka_-
z}\ps  -\fr{\ri\ka_- p_j (z_- - z)}{z_-} t e^{-\ri\ka_- z}(\ps_{xx} + \ps_{yy}).
\]

The transparent boundary conditions \eqref{rtbcp} and \eqref{rtbcm} become
\begin{align}\lb{rtbcp_1}
\left\{
\begin{array}{c}
\pa_z E_1^{+ (1)}(\rh, z_+)-\pa_x E_3^{+ (1)}(\rh, z_+)=\ri\om\mu
H_1^{+ (1)}(\rh, z_+) - f^{+ (1)}_1(\rh),\\[10pt]
\pa_z E_2^{+ (1)}(\rh, z_+)-\pa_y E_3^{+ (1)}(\rh, z_+)=\ri\om\mu
H_2^{+ (1)}(\rh, z_+) - f^{+ (1)}_2(\rh),
\end{array}
\right.
\end{align}
and
\begin{align}\lb{rtbcm_1}
\left\{
\begin{array}{c}
\pa_z E_1^{- (1)}(\rh, z_-)-\pa_x E_3^{- (1)}(\rh, z_-)=-\ri\om\mu
H_1^{- (1)}(\rh, z_-) - f^{- (1)}_1(\rh),\\[10pt]
\pa_z E_2^{- (1)}(\rh, z_-)-\pa_y E_3^{- (1)}(\rh, z_-)=-\ri\om\mu
H_2^{- (1)}(\rh, z_-) - f^{- (1)}_2(\rh),
\end{array}
\right.
\end{align}
where we have from \eqref{s0} that
\begin{align*}
f_j^{+ (1)}(\rh)&=\fr{\ps}{z_+}\pa_z E_j^{+ (0)}(\rh, z_+)=-\fr{\ri\ka_+
p_j}{z_+}\left(e^{-\ri\ka_+ z_+} - re^{\ri\ka_+ z_+}\right)\ps,\\[5pt]
f_j^{- (1)}(\rh)&=\fr{\ps}{z_-}\pa_z E_j^{- (0)}(\rh, z_-)=-\fr{\ri\ka_-
p_j}{z_-}t e^{-\ri\ka_- z_-}\ps.
\end{align*}

The jump conditions \eqref{rce} and \eqref{rch} reduce to
\begin{align*}
\left\{
\begin{array}{c}
E_2^{+ (1)} +\ps_y E_3^{+ (0)} =E_2^{- (1)} +\ps_y E_3^{- (0)},\\[5pt]
E_1^{+ (1)} +\ps_x E_3^{+ (0)} =E_1^{- (1)} +\ps_x E_3^{- (0)},
\end{array}
\right.
\end{align*}
and
\begin{align*}
\left\{
\begin{array}{c}
\pa_z E_1^{+ (1)} -\df{\ps}{z_-}\pa_z E_1^{+ (0)}- \pa_x E_3^{+ (1)}=\pa_z
E_1^{- (1)} - \df{\ps}{z_+}\pa_z E_1^{- (0)}-\pa_x E_3^{- (1)},\\[10pt]
\pa_z E_2^{+ (1)} -\df{\ps}{z_-} \pa_z E_2^{+ (0)} - \pa_y E_3^{+ (1)}=\pa_z
E_2^{- (1)} -\df{\ps}{z_+}\pa_z E_2^{- (0)}- \pa_y E_3^{- (1)},
\end{array}
\right.
\end{align*}
which gives after substitution of \eqref{s0} that
\be\lb{rce_1}
E_2^{+ (1)}  = E_2^{- (1)},\qu E_1^{+ (1)}  = E_1^{- (1)},
\ee
and
\begin{align}\lb{rch_1}
\left\{
\begin{array}{c}
\pa_z E_1^{+ (1)} -\pa_x E_3^{+
(1)} + \df{\ri\ka_+ p_1}{z_-}(1-r)\ps =\pa_z E_1^{- (1)}-\pa_x E_3^{-
(1)} + \df{\ri\ka_-p_1}{z_+}t\ps ,\\[10pt]
\pa_z E_2^{+ (1)} - \pa_y E_3^{+
(1)} +\df{\ri\ka_+ p_2}{z_-}(1-r)\ps =\pa_z E_2^{- (1)} - \pa_y E_3^{-
(1)}+\df{\ri\ka_- p_2}{z_+}t\ps.
\end{array}
\right.
\end{align}

The divergence free condition \eqref{rdf} reduces to
\be\lb{rdf_1}
\pa_x E_1^{\pm (1)}+\pa_y E_2^{\pm (1)}+\pa_z E_3^{\pm
(1)}=g^{\pm (1)}\qu\text{in}~ D_\pm,
\ee
where
\[
g^{\pm (1)}(\rh, z)=\fr{\ps}{z_\pm}\left(\pa_x E_1^{\pm (0)}+\pa_y
E_2^{\pm (0)}\right)+\left(\fr{z_\pm-z}{z_\pm}\right)\left(\ps_x \pa_z
E_1^{\pm (0)}+\ps_y \pa_z E_2^{\pm (0)}\right).
\]
Using \eqref{s0}, we get
\[
 g^{+ (1)}(\rh, z)=-\fr{\ri\ka_+(z_+ -z)}{z_+}\left(e^{-\ri\ka_+ z} -
re^{\ri\ka_+ z}\right)(p_1\ps_x + p_2\ps_y)
\]
and
\[
 g^{- (1)}(\rh, z)=-\fr{\ri\ka_-(z_- -z)}{z_-} t e^{-\ri\ka_- z}(p_1\ps_x +
p_2\ps_y).
\]

Since $\ps(\rh), E_j^{\pm (1)}(\rh, z)$, and $F_j^{\pm (1)}(\rh, z)$ are
periodic functions of $\rh$ with period $\La$, they have the following Fourier
expansions
\begin{align*}
 \ps(\rh)&=\sum_{n\in\mb{Z}^2}\ps_n e^{\ri\al_n\cdot\rh},\\[5pt]
 E_j^{\pm (1)}(\rh, z)&=\sum_{n\in\mb{Z}^2}E_{j
n}^{\pm (1)}(z)e^{\ri\al_n\cdot\rh},\\[5pt]
 F_j^{\pm (1)}(\rh, z)&=\sum_{n\in\mb{Z}^2}F_{j
n}^{\pm (1)}(z)e^{\ri\al_n\cdot\rh}
\end{align*}
where
\begin{align*}
 F_{j n}^{+ (1)}(z)=&\Bigl[\fr{2\ka^2_+ p_j}{z_+}\left(e^{-\ri\ka_+
z}+ re^{\ri\ka_+ z}\right)\\
&+\fr{\ri\ka_+ p_j(z_+ -z)}{z_+}\,(\al_{1 n}^2+\al_{2 n}^2)\left(e^{-\ri\ka_+
z}-re^{\ri\ka_+ z}\right)\Bigr]\ps_n
\end{align*}
and
\[
 F_{j n}^{- (1)}(z)=\Bigl[\fr{2\ka_-^2 p_j}{z_-} t e^{-\ri\ka_-
z}+\fr{\ri\ka_- p_j (z_- - z)}{z_-}(\al_{1 n}^2+\al_{2 n}^2)t e^{-\ri\ka_-
z}\Bigr]\ps_n.
\]

Plugging the above Fourier expansions into \eqref{rtf_1} and using
\eqref{rtbcp_1}--\eqref{rdf_1}, we derive an ordinary differential equation
\be\lb{ode_1f}
\fr{{\rm d}^2 E^{\pm (1)}_{j n}(z)}{{\rm
d}z^2}+ (\ba_n^\pm)^2 E^{\pm (1)}_{j n}(z)=F_{j
n}^{\pm (1)}(z),
\ee
together with the boundary conditions at $z=z_+$:
\begin{align}\lb{rbcp_1f}
\left\{
\begin{array}{l}
{E_{1 n}^{+ (1)}}'-\ri\al_{1 n} E_{3
n}^{+ (1)}=\df{\ri}{\ba^+_n}\left[(\ka_+^2-\al_{2
n}^2)E^{+ (1)}_{1 n}+\al_{1 n}\al_{2 n}E^{+ (1)}_{2 n}\right] - f^{+ (1)}_{1
n},\\[12pt]
{E_{2 n}^{+ (1)}}'-\ri\al_{2 n} E_{3
n}^{+ (1)}=\df{\ri}{\ba^+_n}\left[(\ka_+^2-\al_{1
n}^2)E^{+ (1)}_{2 n}+\al_{1 n}\al_{2 n}E^{+ (1)}_{1 n}\right] - f^{+ (1)}_{2
n}\\[12pt]
{E_{3 n}^{+ (1)}}'+\ri\al_{1 n}E_{1 n}^{+ (1)}+\ri\al_{2 n}E_{2
n}^{+ (1)}=0,
\end{array}
\right.
\end{align}
and the boundary conditions at $z=z_-$:
\begin{align}\lb{rbcm_1f}
\left\{
\begin{array}{l}
{E_{1 n}^{- (1)}}'-\ri\al_{1 n} E_{3
n}^{- (1)}=-\df{\ri}{\ba^-_n}\left[(\ka_-^2-\al_{2
n}^2)E^{- (1)}_{1 n}+\al_{1 n}\al_{2 n}E^{- (1)}_{2 n}\right] - f^{- (1)}_{1
n},\\[12pt]
{E_{2 n}^{- (1)}}'-\ri\al_{2 n} E_{3
n}^{- (1)}=-\df{\ri}{\ba^-_n}\left[(\ka_-^2-\al_{1
n}^2)E^{- (1)}_{2 n}+\al_{1 n}\al_{2 n}E^{- (1)}_{1 n}\right] - f^{-
(1)}_{2 n}\\[12pt]
{E_{3 n}^{- (1)}}'+\ri\al_{1 n}E_{1 n}^{- (1)}+\ri\al_{2 n}E_{2
n}^{- (1)}=0,
\end{array}
\right.
\end{align}
where $f^{\pm (1)}_{j n}$ are the Fourier coefficients of $f^{\pm
(1)}_{j}(\rh)$. Explicitly, we have
\begin{align*}
f_{j n}^{+ (1)}&=-\fr{\ri\ka_+
p_j}{z_+}\left(e^{-\ri\ka_+ z_+} - re^{\ri\ka_+ z_+}\right)\ps_n,\\[5pt]
f_{j n}^{- (1)}&=-\fr{\ri\ka_- p_j}{z_-}t e^{-\ri\ka_- z_-}\ps_n.
\end{align*}

Using the identity $\ka_+ (1-r)=\ka_- t$, we may reduce the jump conditions
\eqref{rce_1} and \eqref{rch_1} to
\begin{equation}\lb{rce_1f}
E_{2 n}^{+ (1)} =E_{2 n}^{- (1)},\qu E_{1 n}^{+ (1)} =E_{1 n}^{- (1)},
\end{equation}
and
\begin{align}\lb{rch_1f}
\left\{
\begin{array}{c}
 {E_{1 n}^{+ (1)}}' -\ri\al_{1 n} E_{3 n}^{+
(1)} ={E_{1 n}^{- (1)}}' -\ri\al_{1 n}
E_{3 n}^{- (1)}+\ri\ka_-
t p_1\left(z_+^{-1}-z_-^{-1}\right)\ps_n,\\[5pt]
{E_{2 n}^{+ (1)}}' - \ri\al_{2 n} E_{3 n}^{+ (1)}={E_{2 n}^{- (1)}}' - \ri\al_{2
n}E_{3 n}^{- (1)} +\ri\ka_-
t p_2\left(z_+^{-1}-z_-^{-1}\right)\ps_n
\end{array}
\right.
\end{align}
Based on the same identity $\ka_+ (1-r)=\ka_- t$, we may obtain two more
conditions at $z=0$ from \eqref{rdf_1}:
\begin{align}\lb{rcd_1f}
\left\{
\begin{array}{l}
 {E_{3 n}^{+ (1)}}'+\ri\al_{1 n}E_{1 n}^{+ (1)}+\ri\al_{2 n}E_{2
n}^{+ (1)}=\ka_- t (\al_{1 n}p_1 + \al_{2 n}p_2) \ps_n,\\[5pt]
{E_{3 n}^{- (1)}}'+\ri\al_{1 n}E_{1 n}^{- (1)}+\ri\al_{2 n}E_{2
n}^{- (1)}=\ka_- t (\al_{1n}p_1 + \al_{2 n}p_2) \ps_n.
\end{array}
\right.
\end{align}

It follows from \eqref{ode_1f} that the general solutions of $E_{j
n}^{\pm (1)}$ consist of the general solution for the corresponding
homogeneous equation and a particular solution for the non-homogeneous equation:
\be\lb{s1p}
E_{j n}^{+ (1)}(z)=A^+_{j n} e^{\ri\ba^+_n z}+B^+_{j n} e^{-\ri\ba^+_n
z}-\df{\ri\ka_+ p_j}{z_+}(z_+ -z)\left(e^{-\ri\ka_+ z}- r
e^{\ri\ka_+ z}\right)\ps_n
\ee
and
\be\lb{s1m}
E_{j n}^{- (1)}(z)=A^-_{j n} e^{\ri\ba^-_n z}+B^-_{j n} e^{-\ri\ba^-_n
z}-\df{\ri\ka_- p_j}{z_-}(z_- -z) t e^{-\ri\ka_- z} \ps_n
\ee

Plugging \eqref{s1p} and \eqref{s1m} into \eqref{rbcp_1f} and \eqref{rbcm_1f},
and using the identity $\ka_\pm^2=(\ba^\pm_n)^2+\al_{1 n}^2+\al_{2 n}^2$, we
obtain
\begin{align}\lb{ABp}
\left\{
\begin{array}{l}
\al_{1 n}^2 A^+_{1 n}+[2(\ba^+_n)^2+\al_{1 n}^2]e^{-2\ri\ba^+_n z_+}B^+_{1
n}+\al_{1 n}\al_{2 n}A^+_{2 n}\\[5pt]
\hs{1cm}+\al_{1 n}\al_{2 n}e^{-2\ri\ba^+_n
z_+}B^+_{2 n}=-\al_{1 n}\ba^+_n A^+_{3 n}-\al_{1 n}\ba^+_n e^{-2\ri\ba^+_n
z_+}B^+_{3
n},\\[10pt]
\al_{2 n}^2 A^+_{2 n}+[2(\ba^+_n)^2+\al_{2 n}^2]e^{-2\ri\ba^+_n z_+}B^+_{2
n}+\al_{1 n}\al_{2 n}A^+_{1 n}\\[5pt]
\hs{1cm}+\al_{1 n}\al_{2 n}e^{-2\ri\ba^+_n z_+}B^+_{1
n}=-\al_{2 n}\ba_n A^+_{3 n}-\al_{2 n}\ba^+_n e^{-2\ri\ba^+_n z_+}B^+_{3
n},\\[10pt]
\al_{1 n}A^+_{1 n}+\al_{1 n}e^{-2\ri\ba^+_n z_+}B^+_{1 n}+\al_{2 n}A^+_{2
n}+\al_{2 n}e^{-2\ri\ba^+_n z_+}B^+_{2 n}\\[5pt]
\hs{1cm}=-\ba^+_n A^+_{3 n}+\ba^+_n e^{-2\ri\ba^+_n z_+} B^+_{3  n}.
\end{array}
\right.
\end{align}
and
\begin{align}\lb{ABm}
\left\{
\begin{array}{l}
[2(\ba^-_n)^2+\al_{1 n}^2] A^-_{1 n}+ \al_{1 n}^2 e^{-2\ri\ba^-_n z_-}B^-_{1
n}+\al_{1 n}\al_{2 n}A^-_{2 n}\\[5pt]
\hs{1cm}+\al_{1 n}\al_{2 n}e^{-2\ri\ba^-_n z_-}B^-_{2
n}=\al_{1 n}\ba^-_n A^-_{3 n}+\al_{1 n}\ba^-_n e^{-2\ri\ba^-_n z_-}B^-_{3
n},\\[10pt]
[2(\ba^-_n)^2+\al_{2 n}^2]  A^-_{2 n}+\al_{2 n}^2 e^{-2\ri\ba^-_n z_-}B^-_{2
n}+\al_{1 n}\al_{2 n}A^-_{1 n}\\[5pt]
\hs{1cm}+\al_{1 n}\al_{2 n}e^{-2\ri\ba^-_n z_-}B^-_{1
n}=\al_{2 n}\ba^-_n A^-_{3 n}+\al_{2 n}\ba^-_n e^{-2\ri\ba^-_n z_-}B^-_{3
n},\\[10pt]
\al_{1 n}A^-_{1 n}+\al_{1 n}e^{-2\ri\ba^-_n z_-}B^-_{1 n}+\al_{2 n}A^-_{2
n}+\al_{2 n}e^{-2\ri\ba^-_n z_-}B^-_{2 n}\\[5pt]
\hs{1cm}=-\ba^-_n A^-_{3 n}+\ba^-_n e^{-2\ri\ba^-_n z_-} B^-_{3  n}.
\end{array}
\right.
\end{align}
Multiplying individually $\al_{1 n}$ and $\al_{2 n}$ on both sides of the third
equation in \eqref{ABp} and \eqref{ABm}, and subtracting them from the first and
second equation, respectively, we get
\be\lb{AmBp}
B^+_{1 n}=-\df{\al_{1 n}}{\ba^+_n} B^+_{3 n},\qu B^+_{2
n}=-\df{\al_{2 n}}{\ba^+_n} B^+_{3 n},\qu A^-_{1 n}=\df{\al_{1 n}}{\ba^-_n}
A^-_{3 n}, \qu A^-_{2 n}=\df{\al_{2 n}}{\ba^-_n} A^-_{3 n}.
\ee
Substituting \eqref{AmBp} into the third equations in \eqref{ABp} and
\eqref{ABm} yields
\begin{align}\lb{AB_e1}
\left\{
\begin{array}{c}
\al_{1 n}A_{1 n}^+ + \al_{2 n}A_{2 n}^+=\df{\ka_+^2}{\ba_n^+}e^{-2\ri\ba_n^+
z_+} B_{3 n}^+ - \ba_n^+ A_{3 n}^+,\\[5pt]
\al_{1 n}B_{1 n}^- + \al_{2 n} B_{2 n}^- =\ba_n^- B_{3 n}^-
-\df{\ka_-^2}{\ba_n^-}e^{2\ri\ba_n^- z_-} A_{3 n}^-.
\end{array}
\right.
\end{align}
Substituting \eqref{AmBp}, \eqref{s1p}, and \eqref{s1m} into \eqref{rcd_1f}, we
get
\begin{align}\lb{AB_e2}
\left\{
\begin{array}{c}
\al_{1 n}A_{1 n}^+ + \al_{2 n}A_{2 n}^+=\df{\ka_+^2}{\ba_n^+} B_{3 n}^+ -
\ba_n^+ A_{3 n}^+,\\[5pt]
\al_{1 n}B_{1 n}^- + \al_{2 n} B_{2 n}^- =\ba_n^- B_{3 n}^-
-\df{\ka_-^2}{\ba_n^-} A_{3 n}^-.
\end{array}
\right.
\end{align}
Combining \eqref{AmBp}--\eqref{AB_e2} gives
\be\lb{AB1}
B_{1 n}^+=B_{2 n}^+=B_{3 n}^+=0\qu\text{and}\qu A_{1 n}^-=A_{2 n}^-=A_{3 n}^-=0.
\ee
Plugging \eqref{AB1}, \eqref{s1p}, and \eqref{s1m} into \eqref{rce_1f} and
\eqref{rch_1f}, we obtain
\begin{align}\lb{AB_e3}
\left\{
\begin{array}{c}
A_{1 n}^+ - B_{1 n}^- = 0,\\[5pt]
\ba_n^+ A_{1 n}^+ +\ba_n^- B_{1 n}^- = \al_{1 n} (A_{3 n}^+ - B_{3
n}^-)-2\ri\ka_+(\ka_+ - \ka_-) p_1 \ps_n,
\end{array}
\right.
\end{align}
and
\begin{align}\lb{AB_e4}
\left\{
\begin{array}{c}
A_{2 n}^+ - B_{2 n}^- = 0,\\[5pt]
\ba_n^+ A_{2 n}^+ +\ba_n^- B_{2 n}^- = \al_{2 n} (A_{3 n}^+ - B_{3
n}^-)-2\ri\ka_+(\ka_+ - \ka_-) p_2 \ps_n.
\end{array}
\right.
\end{align}
Upon solving \eqref{AB_e3} and \eqref{AB_e4}, we have
\begin{align}\lb{AB2}
\left\{
\begin{array}{c}
A_{1 n}^+ = B_{1 n}^-= (\ba_n^+ + \ba_n^-)^{-1} \left[\al_{1 n} (A_{3 n}^+ -
B_{3 n}^-)-2\ri\ka_+(\ka_+ - \ka_-) p_1 \ps_n\right],\\[10pt]
A_{2 n}^+ = B_{2 n}^- =(\ba_n^+ + \ba_n^-)^{-1}\left[\al_{2 n} (A_{3 n}^+ - B_{3
n}^-)-2\ri\ka_+(\ka_+ - \ka_-) p_2 \ps_n\right].
\end{array}
\right.
\end{align}
Substituting \eqref{AB2} into \eqref{AB_e2} and noting \eqref{AB1}, we may
derive after tedious calculations that
\begin{align}\lb{AB3}
\left\{
\begin{array}{l}
A_{3 n}^+ = \df{2\ri\ba_n^-\ka_+ (\ka_+ - \ka_-)(p_1\al_{1
n}+p_2\al_{2 n})}{(\ba_n^+ + \ba_n^-)(\al_{1 n}^2 + \al_{2 n}^2 + \ba_n^+
\ba_n^-)}\ps_n,\\[10pt]
B_{3 n}^- = -\df{2\ri\ba_n^+ \ka_+ (\ka_+ - \ka_-)(p_1\al_{1
n}+p_2\al_{2 n})}{(\ba_n^+ + \ba_n^-)(\al_{1 n}^2 + \al_{2 n}^2 + \ba_n^+
\ba_n^-)}\ps_n.
\end{array}
\right.
\end{align}
Plugging \eqref{AB3} into \eqref{AB2} yields
\be\lb{AB4}
A_{1 n}^+=B_{1 n}^-=C_{1 n}\ps_n, \qu A_{2 n}^+=B_{2 n}^-=C_{2 n}\ps_n,
\ee
where
\begin{align*}
\left\{
\begin{array}{l}
C_{1 n}=\df{2\ri\ka_+ (\ka_+ - \ka_-)}{(\ba_n^+ +
\ba_n^-)}\left[\df{\al_{1 n}(p_1\al_{1 n}+ p_2\al_{2 n})}{(\al_{1 n}^2 + \al_{2
n}^2 + \ba_n^+ \ba_n^-)} - p_1\right],\\[10pt]
C_{2 n}=\df{2\ri\ka_+ (\ka_+ - \ka_-)}{(\ba_n^+ +
\ba_n^-)}\left[\df{\al_{2 n}(p_1\al_{1 n}+ p_2\al_{2 n})}{(\al_{1 n}^2 + \al_{2
n}^2 + \ba_n^+ \ba_n^-)} - p_2\right].
\end{array}
\right.
\end{align*}
Substituting \eqref{AB1} and \eqref{AB4} into \eqref{s1p} and \eqref{s1m}, and
evaluating at
$z_+$ and $z_-$, respectively, we obtain
\be\lb{s1}
E_{j n}^{+ (1)}(z_+)=C_{j n} e^{\ri\ba_n^+ z_+}\ps_n,\qu E_{j n}^{-
(1)}(z_-)=C_{j n} e^{-\ri\ba_n^- z_-}\ps_n.
\ee

\section{Reconstruction formula}

In this section, we present an explicit reconstruction formula for the inverse
grating surface problem by using the scattering data.

Assume that the noisy data takes the form
\[
E_j^{\pm \ga}(\rh, z_\pm)=E^\pm_j(\rh, z_\pm)+\mc{O}(\ga),
\]
where $E^\pm_j(\rh, z_\pm), j = 1, 2$ are the exact data and $\ga$ is the noise
level.

Evaluating the power series \eqref{ps} at $z=z_\pm$ and replacing $E^\pm_j(\rh,
z_\pm)$ with the noisy data $E_j^{\pm \ga}(\rh, z_\pm)$, we have
\be\lb{wd}
E^{\pm \ga}_j(\rh, z_\pm)=E^{\pm (0)}_j(\rh, z_\pm)+\de E^{\pm (1)}_j(\rh,
z_\pm)+\mc{O}(\de^2)+\mc{O}(\ga).
\ee
Rearranging \eqref{wd}, and dropping $\mc{O}(\de^2)$ and $\mc{O}(\ga)$ yield
\be\lb{rw}
\de E_j^{\pm (1)}(\rh, z_\pm)=E^{\pm \ga}_j(\rh, z_\pm)-E_j^{\pm (0)}(\rh,
z_\pm)
\ee
which is the linearization of the nonlinear inverse problem and enables us to
find an explicit reconstruction formula for the linearized inverse problem.

Noting $\ph=\de\ps$ and thus $\ph_n=\de\ps_n$, where $\ph_n$ is the Fourier
coefficient of $\ph$. Plugging \eqref{s1} into \eqref{rw}, we may
deduce that
\be\lb{rf}
\ph_n=C_{j n}^{-1}\left[E_{jn}^{\pm
\ga}(z_\pm)-E_{j n}^{\pm (0)}(z_\pm)\right]e^{\mp\ri\ba^\pm_n
z_\pm},
\ee
where $E_{jn}^{\pm \ga}(z_\pm)$ is the Fourier coefficient of the noisy data
$E^{\pm \ga}_j(\rh, z_\pm)$ and $E_{jn}^{\pm (0)}(z_\pm)$ is the Fourier
coefficient of $E_j^{\pm (0)}(\rh, z_\pm)$ given as
\be\lb{s0pm}
E_{j n}^{+(0)}(z_+)=p_j(e^{-\ri\ka_+ z_+}+ r e^{\ri\ka_+ z_+})\de_{0
n}\qu\text{and}\qu E_{j n}^{-(0)}(z_-)=p_j t e^{-\ri\ka_- z_-} \de_{0
n}.
\ee
Here $\de_{0n}$ the Kronecker's delta function.

It follows from \eqref{rf} and the definitions of $\ba^\pm_n$ in \eqref{bap}, 
\eqref{bam} that it is well-posed to reconstruct those Fourier coefficients
$\ph_n$ with $|\al_n|<\ka_\pm$, since the small variations of the measured data
will not be amplified and lead to large errors in the reconstruction, but the
resolution of the reconstructed function $f$ is restricted by the given
wavenumber $\ka_\pm$. In contrast, it is severely ill-posed to reconstruct those
Fourier coefficients $\ph_n$ with $|\al_n|>\ka_\pm$, since the small variations
in the data will be exponentially enlarged and lead to huge errors in the
reconstruction, but they contribute to the super resolution of the reconstructed
function $\ph$.

To obtain a stable and super-resolved reconstruction, we adopt
a regularization to suppress the exponential growth of the reconstruction
errors. Besides, we may use as small $|z_\pm|$ as possible, i.e., measure the
data at the distance which is as close as possible to the grating surface which
is exactly the idea of near-field optics.

We consider the spectral cut-off regularization. Define the signal-to-noise
ratio (SNR) by
\[
{\rm SNR}=\min\{\de^{-2}, ~\ga^{-1}\}.
\]
For fixed $z_\pm$, the cut-off frequency $\om_\pm$ is chosen in such a way that
\[
e^{ |z_\pm| (\om_\pm^2-\ka_\pm^2)^{1/2}}={\rm SNR},
\]
which implies that the spatial frequency will be cut-off for those below the
noise level. More explicitly, we have
\[
\fr{\om_\pm}{\ka_\pm}=\left[1+\left(\fr{\log{\rm SNR}}{\ka_\pm
|z_\pm|}\right)^2\right]^{1/2},
\]
which indicates $\om_\pm>\ka_\pm$ as long as ${\rm SNR}>0$ and super resolution
may be achieved.

Taking into account the frequency cut-off, we may have a regularized
reconstruction formulation for \eqref{rf}:
\[
\ph_n=C_{j n}^{-1}\left[E_{j n}^{\pm
\ga}(z_\pm) - E_{j n}^{\pm (0)}(z_\pm)\right]e^{\mp\ri\ba^\pm_n
z_\pm}\,\chi^\pm_n,
\]
where the characteristic function
\begin{align*}
\chi^\pm_n=\left\{
\begin{array}{lll}
1 &\qu\text{for} & |\al_n|\le\om_\pm,\\[5pt]
0 &\qu\text{for} & |\al_n|>\om_\pm.
\end{array}
\right.
\end{align*}
Once $\ph_n$ are computed, the grating surface function can be approximated by
\begin{align}\lb{rphn}
\ph(\rh)&\approx\sum_{n\in\mb{Z}}\ph_n
e^{\ri\al_n\cdot\rh}=\sum_{|\al_n|\le\om_\pm}
C_{j n}^{-1}\left[E_{jn}^{\pm \ga}(z_\pm)-E_{j n}^{\pm
(0)}(z_\pm)\right]e^{\ri(\al_n\cdot\rh\mp\ba^\pm_n
z_\pm)}\nt\\[5pt]
&=\sum_{|\al_n|\le\om_\pm} C_{j n}^{-1}
E_{jn}^{\pm \ga}(z_\pm)e^{\ri(\al_n\cdot\rh\mp\ba^\pm_n
z_\pm)}-\sum_{|\al_n|\le\om_\pm} C_{j n}^{-1}
E_{jn}^{\pm (0)}(z_\pm)e^{\ri(\al_n\cdot\rh\mp\ba^\pm_n
z_\pm)}.
\end{align}
Substituting \eqref{s0pm} into \eqref{rphn}, we obtain an reconstructed grating
surface function
\[
\ph(\rh)\approx\sum_{|\al_n|\le\om_\pm} C_{j n}^{-1}
E_{jn}^{\pm \ga}(z_\pm)e^{\ri(\al_n\cdot\rh\mp\ba^\pm_n
z_\pm)}- C_{j 0}^{-1} \left(r + e^{-2\ri\ka_+ z_+}\right)p_j
\]
from the reflection configuration or
\[
\ph(\rh)\approx \sum_{|\al_n|\le\om_\pm} C_{j n}^{-1}
E_{jn}^{\pm \ga}(z_\pm)e^{\ri(\al_n\cdot\rh\mp\ba^\pm_n
z_\pm)}- C_{j 0}^{-1} t p_j
\]
from the transmission configuration.

Hence, only two fast Fourier transforms are needed to reconstruct the grating
surface function: one is done for the data to obtain $E^{\pm \ga}_{j n}(z_\pm)$
and another is done to obtain the approximated function $\ph$.

\section{Numerical experiment}

\begin{figure}
\center
\includegraphics[width=0.4\textwidth]{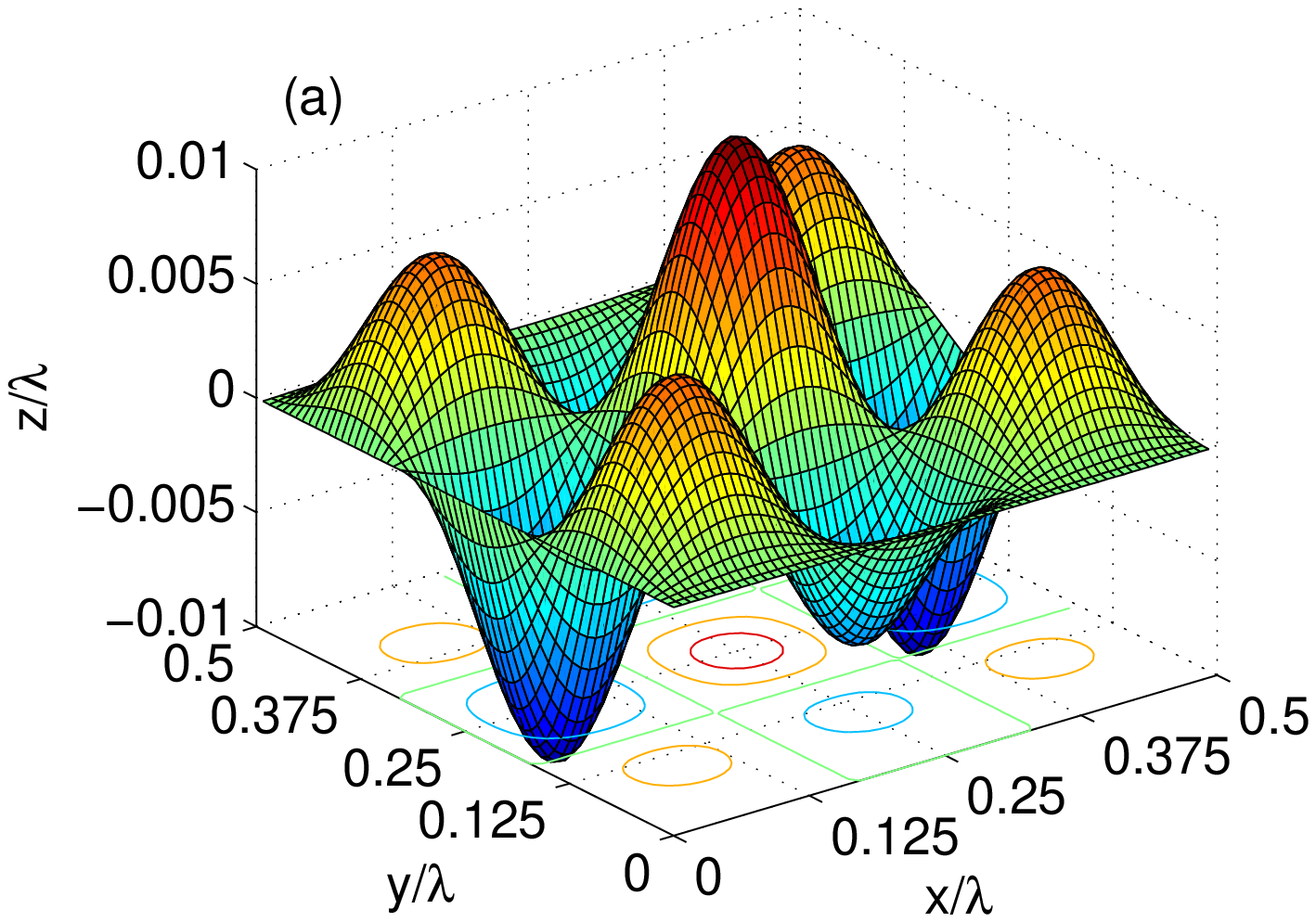}
\includegraphics[width=0.4\textwidth]{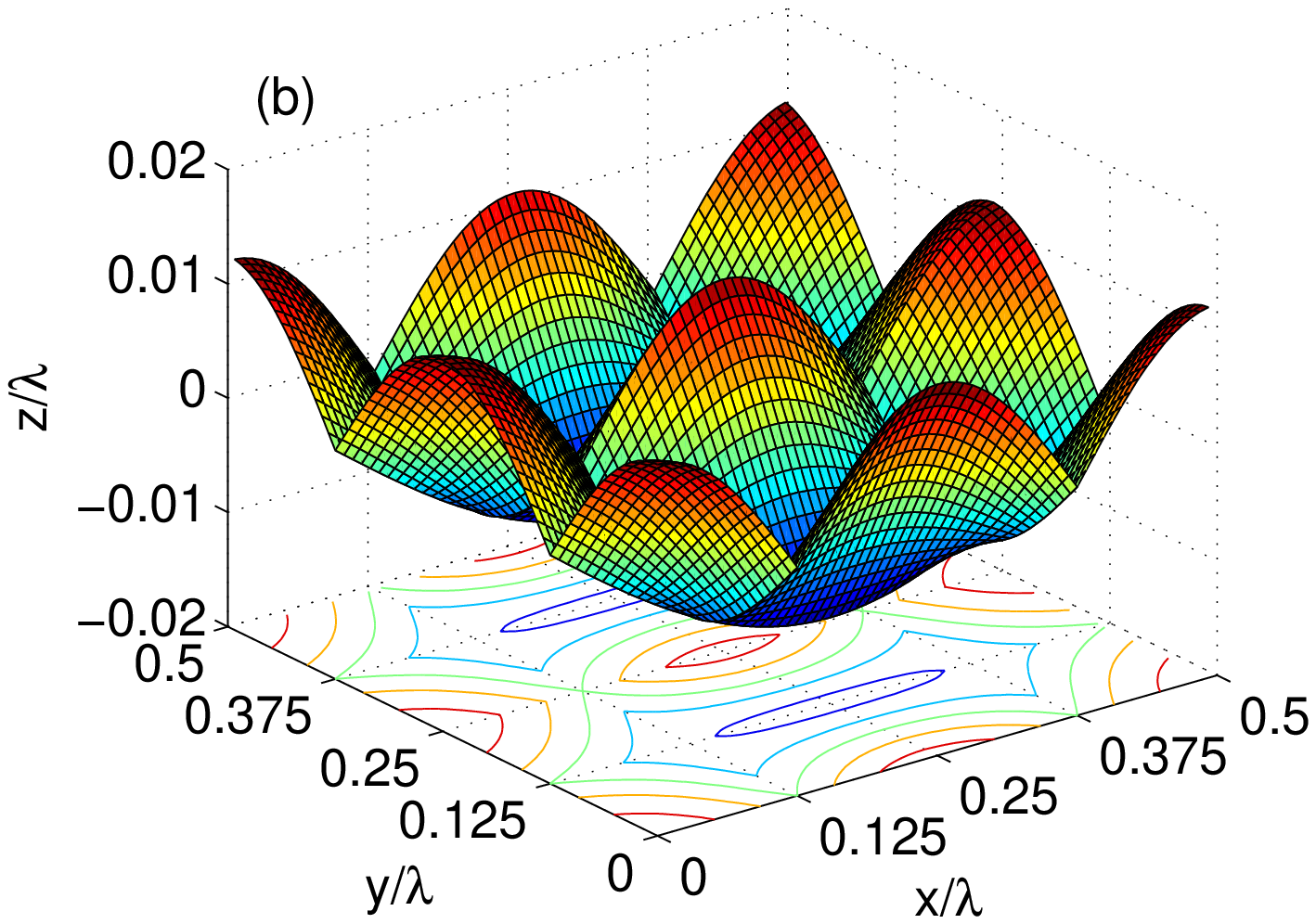}
\caption{The exact grating profile $\ps$. (a) Example 1: smooth grating profile
with finitely many Fourier modes; (b) Example 2: non-smooth grating profile with
infinitely many Fourier modes.}
\label{fig:pro}
\end{figure}

In this section, we discuss the algorithmic implementation for the direct and
inverse problems and present two numerical examples to illustrate the
effectiveness of the proposed method. As is shown in Fig. \ref{fig:pro},
two types of grating profiles are considered. One is a smooth function with
finitely many Fourier modes and another is a non-smooth function with infinitely
many Fourier modes. Although the method requires $\psi\in C^2(\mathbb{R}^2)$, it
is applicable to non-smooth functions numerically.

The second-order N\'{e}d\'{e}lec edge element is adopted to solve the
direct problem and obtain the synthetic scattering data. Uniaxial perfect
matched layer (PML) boundary condition is imposed on the $z$ direction to
truncated the domain. An adaptive mesh refinement technique \cite{B-L-W-1} is
used to achieve the solution with a specified accuracy in an optimal fashion.
Our implementation is based on parallel hierarchical grid (PHG) \cite{phg},
which is a toolbox for developing parallel adaptive finite element programs on
unstructured tetrahedral meshes. To have a tetrahedral mesh with biperiodic
boundary points, we generate an uniform hexahedral mesh and then divide each
hexahedron into six tetrahedrons. The linear system resulted from
finite element discretization is solve by the multifrontal massively parallel
sparse direct solver \cite{mumps1}.

In the following two examples, the incident wave is taken as ${\bf E}^{\rm inc}
= (1, 0, 0)e^{-\ri\ka_+ z}$, i.e., $p_1 = 1$ and $p_2 = p_3 = 0$, and only the
first component of the electric field, $E^+_1(\rho,h)$, needs to be
measured. The wavenumber in $\Omega^-_S$ is $\ka_-=1.6\pi$.
The wavenumber in $\Omega^+_S$ is $\ka_+=\pi$, which corresponds to the
wavelength $\lambda=2$. Define by $R$ the unit rectangular domain, i.e., $R =
[0, 0.5\lambda]\times[0, 0.5\lambda]$. The computational domain is
$R\times[-0.3\lambda,0.3\lambda]$ with the PML region
$(R\times[-0.3\lambda,-0.15\lambda])\cup (R\times[0.15\lambda,0.3\lambda])$.
The scattering data $E^+_1(\rho,h)$ is obtained by interpolation into the
uniform $256\times 256$ grid points on the measurement plane $z = h$. In all
the figures, the plots are rescaled with respect to the wavelength $\lambda$ to
clearly show the relative size. The results are plotted on $64\times 64$
grid points instead of $256\times 256$ grid points in order to reduce
the display sizes. To test the stability of the method, a random
noise is added to the scattering data, i.e., the scattering data takes the form
\[
E^{+\ga}_1(\rho, h)=E^+_1(\rho,h)(1+\ga\,{\rm rand}),
\]
where rand stands for uniformly distributed random numbers in $[-1, 1]$ and
$\ga$ is the noise level parameter. The relative $L^2(R)$ error is
defined by
\[
e=\frac{\|\phi-\phi_{\ga,\de}\|_{0,R}}{\|\phi\|_{0,R}},
\]
where $\phi$ is the exact surface function and $\phi_{\ga,\de}$ is the
reconstructed surface function.

{\em Example 1}. This example illustrates the reconstruction results of a smooth
grating profile with finitely many Fourier modes, as seen in Figure
\ref{fig:pro}(a). The exact grating surface function is given by
$\phi(\rho)=\de\ps(\rho)$, where the grating profile function
\[
\ps(x,y)=0.5\sin(3\pi x)(\cos(2\pi y)-\cos(4\pi y)).
\]
First, consider the surface deviation parameter $\de$. The measurement is taken
at $h = 0.1\lambda$ and no additional random noise is added to the scattering
data, i.e., $\ga = 0$. This test is to investigate the influence of surface
deformation parameter on the reconstructions. In \eqref{rw}, higher order terms
of $\de$ are dropped in the power series to linearize the inverse problem and to
obtain the explicit reconstruction formulas. As expected, the smaller the
surface deformation $\de$ is, the more accurate is the approximation of the
linearized model to the original nonlinear model problem. Table \ref{tab:ex1}
shows the relative $L^2(R)$ error of the reconstructions with three
different surface deformation parameters $\de = 0.05\lambda, 0.025\lambda,
0.0125\lambda$ for a fixed measurement distance $h = 0.1\lambda$. It is clear
to note that the error decreases from $45.3\%$ to $15.6\%$ as $\de$ decreases
from  $0.05\lambda$ to $0.0125\lambda$.

\begin{table}[http!]
\caption{Example 1: Relative error of the reconstructions by using different
$\de$ with $h=0.1\lambda$
and $\ga = 0.0$.}\label{tab:ex1}
\begin{center}{
\begin{tabular}{llll}
\hline
\hline
$\de$   &  $0.05\lambda$         &  $0.025\lambda$         &  $0.0125\lambda$   
  \\
\hline
$e$          &  $4.53\times 10^{-1}$  &  $2.49\times 10^{-1}$  &  $1.56\times
10^{-1}$  \\
\hline
\hline
\end{tabular}}
\end{center}
\end{table}

\begin{figure}
\center
\includegraphics[width=0.4\textwidth]{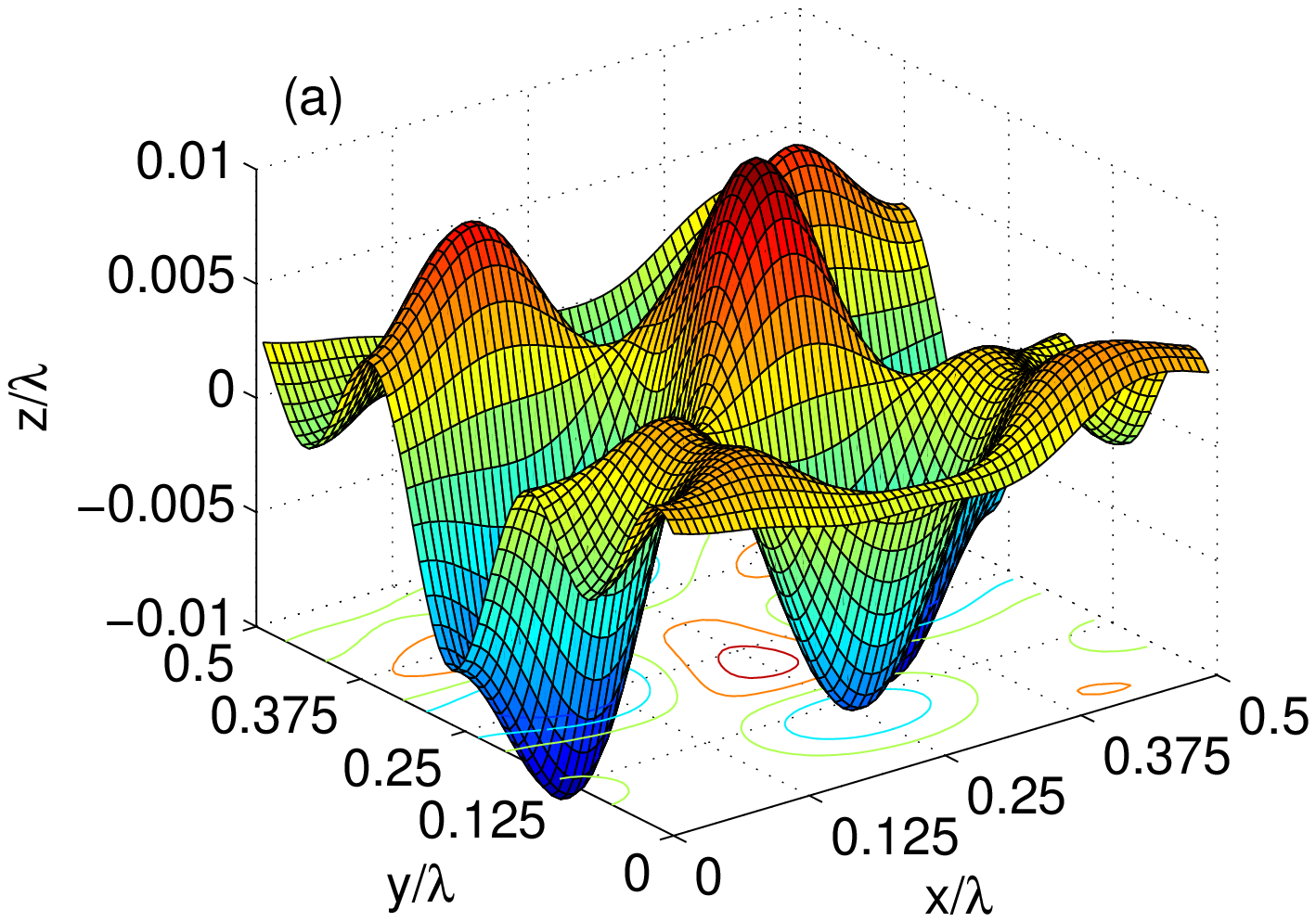}
\includegraphics[width=0.4\textwidth]{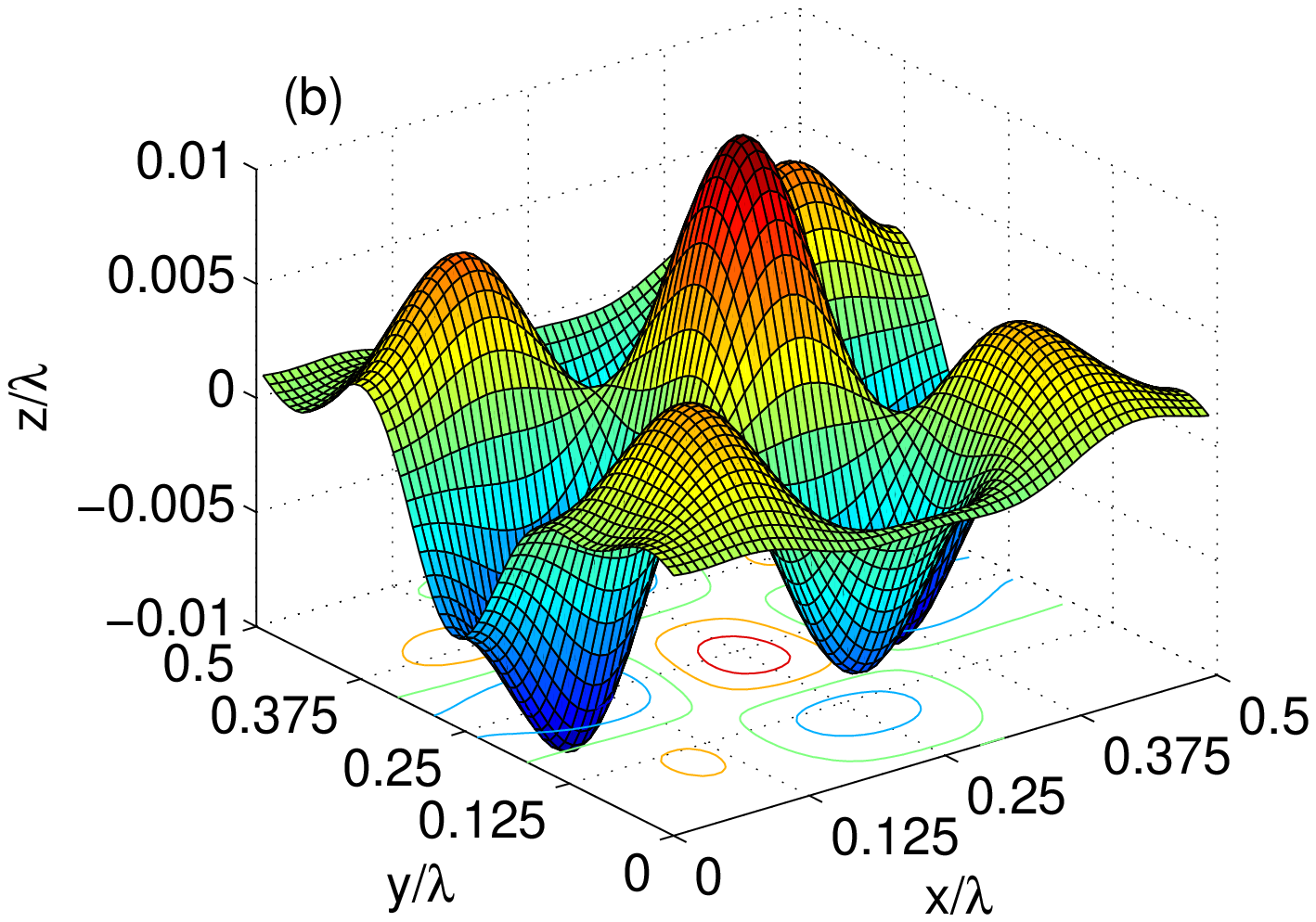}
\includegraphics[width=0.4\textwidth]{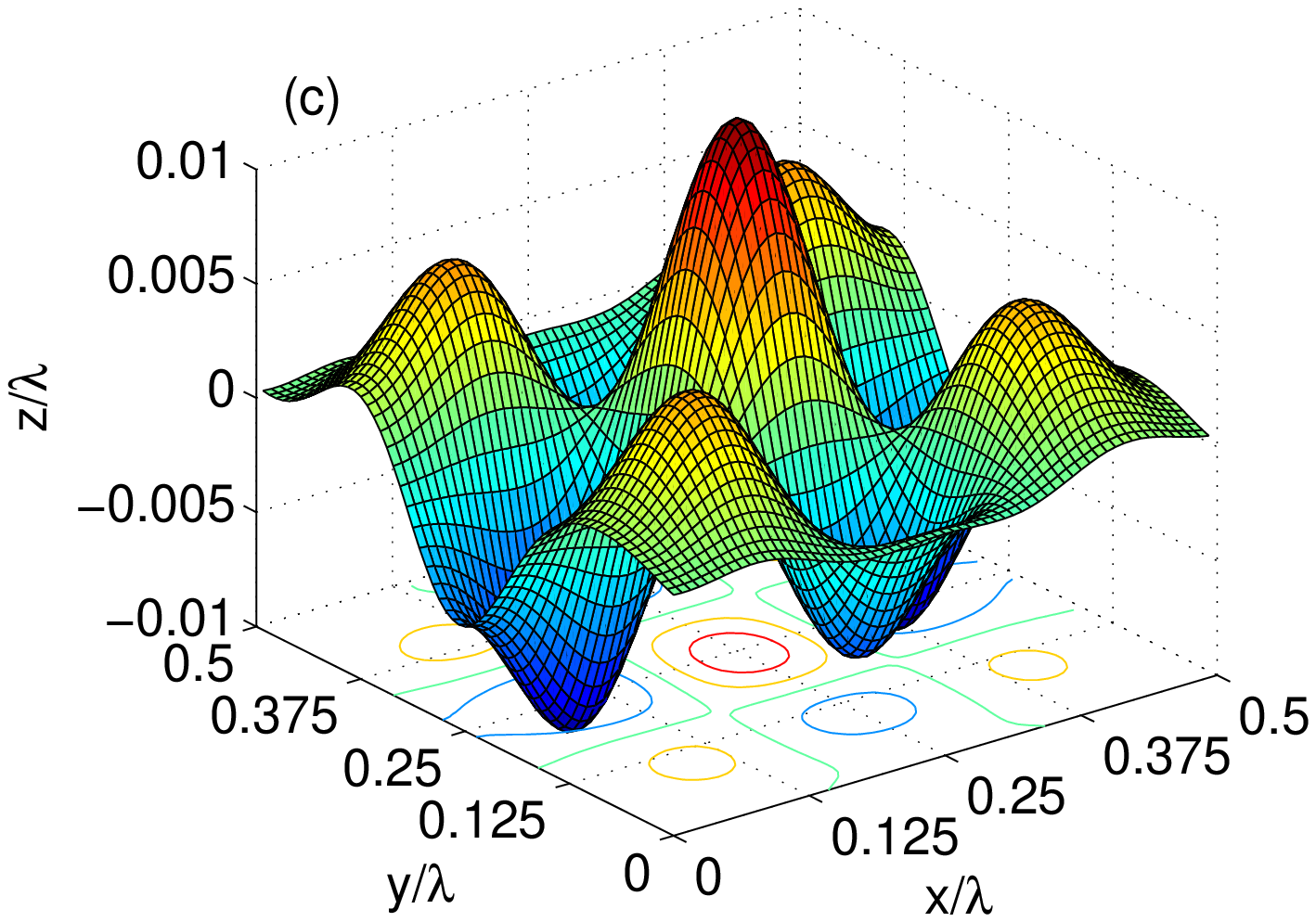}
\includegraphics[width=0.4\textwidth]{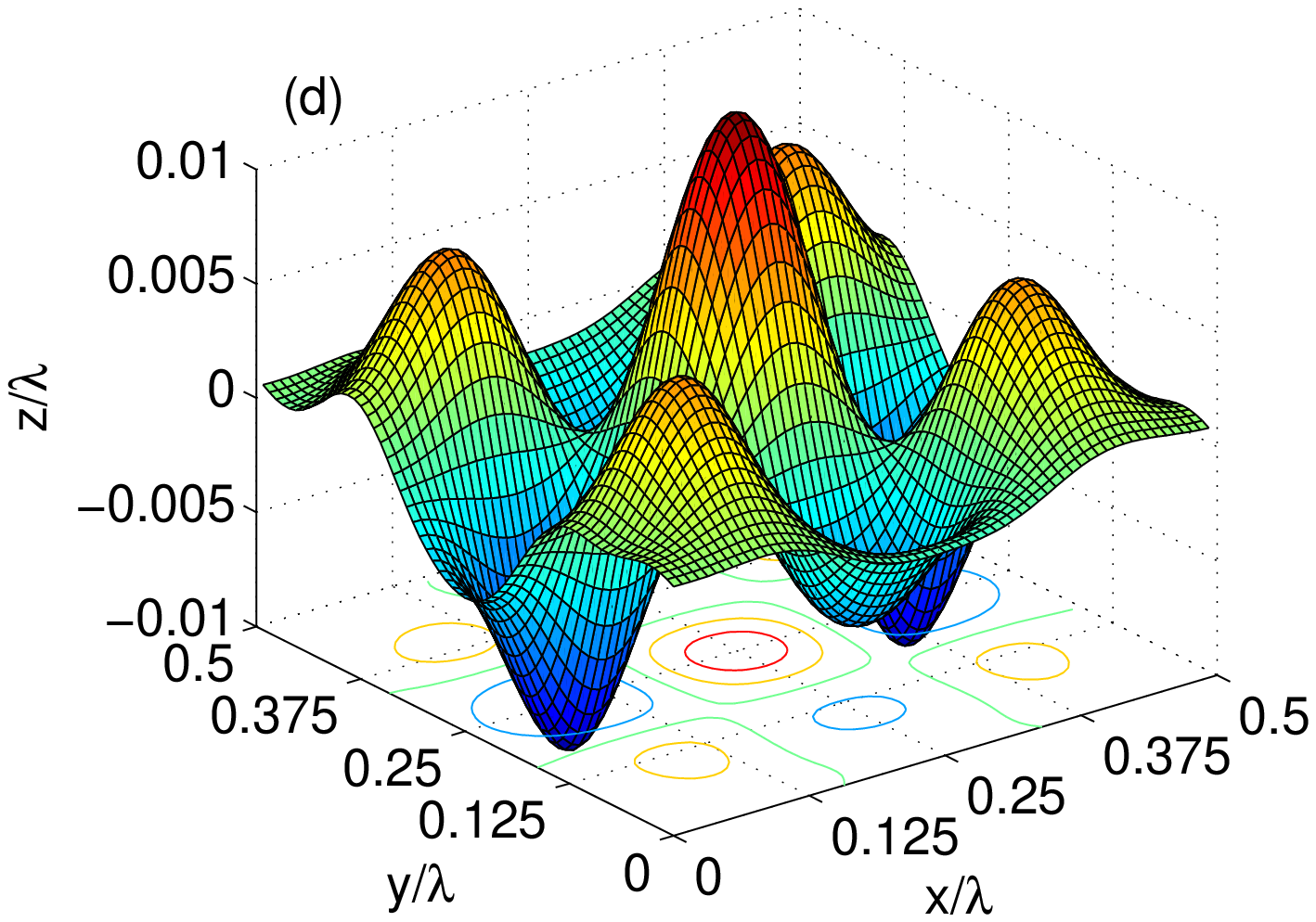}
\caption{Example 1: Reconstructed grating surfaces by using different $h$ with
$\delta = 0.0125\lambda$ and $\ga = 1\%$. (a) $h = 0.1\lambda$; (b) $h =
0.075\lambda$; (c) $h = 0.05\lambda$; (d) $h = 0.025\lambda$.}
\label{fig:ex1}
\end{figure}

Next is to consider the noise level $\ga$ and the measurement distance $h$. In
practice, the scattering data always contains a certain amount of noise. To test
the stability and super resolving capability of the method, we add $1\%$ and
$5\%$ random noises to the scattering data. Table \ref{tab:ex1-1} and
\ref{tab:ex1-2} report the relative $L^2(R)$
error of the reconstructions with four different measurement distances $h
=0.1\lambda, 0.075\lambda, 0.05\lambda, 0.025\lambda$ for a fixed $\de
=0.0125\lambda$. Comparing the results for the same $\de =0.0125\lambda$
and $h = 0.1\lambda$ in Tables \ref{tab:ex1} and \ref{tab:ex1-2}, we can see
that the relative error increases dramatically from $15.6\%$ by using noise free
data to $83.8\%$ by using $5\%$ noise data. The reason is that a smaller
cut-off should be chosen to suppress the exponentially increasing noise in the
data and thus the Fourier modes of the exact grating surface function can not be
recovered for those higher than the cut-off frequency, which leads to a large
error and poor resolution in the reconstruction. A smaller measurement distance
is desirable in order to have a large cut-off frequency, which enhances the
resolution and reduces the error. As can be seen in Table \ref{tab:ex1-1}, the
reconstruction error decreases from $56.7\%$ by using $h = 0.1\lambda$ to as low
as $16.7\%$ by using $h = 0.025\lambda$ for $1\%$ noise data. Similarly, in
Table \ref{tab:ex1-2}, the reconstruction error decreases from $83.8\%$ by using
$h = 0.1\lambda$ to as low as $29.5\%$ by using $h = 0.025\lambda$ even for
$5\%$ noise data. Figure \ref{fig:ex1} plots the reconstructed surfaces by using
$h =0.1\lambda, 0.075\lambda, 0.05\lambda, 0.025\lambda$. Comparing the exact
surface profile in Fig. \ref{fig:pro}(a) and the reconstructed surface in Fig.
\ref{fig:ex1}(d), we can see that the reconstruction is almost perfect and the
difference is little by carefully checking the contour plots.

\begin{table}[http!]
\caption{Example 1: Relative error of the reconstructions by using different $h$
with $\de=0.0125\lambda$
and $\ga = 1\%$.}\label{tab:ex1-1}
\begin{center}{
\begin{tabular}{lllll}
\hline
\hline
$h$    &  $0.1\lambda$     &  $0.075\lambda$         &  $0.05\lambda$         & 
$0.025\lambda$     \\
\hline
$e$    &  $5.67\times 10^{-1}$   &  $2.95\times 10^{-1}$    &  $2.08\times
10^{-1}$  &  $1.67\times 10^{-1}$   \\
\hline
\hline
\end{tabular}}
\end{center}
\end{table}

\begin{table}[http!]
\caption{Example 1: Relative error of the reconstructions by using different $h$
with $\de=0.0125\lambda$
and $\ga = 5\%$.}\label{tab:ex1-2}
\begin{center}{
\begin{tabular}{lllll}
\hline
\hline
$h$    &  $0.1\lambda$    &  $0.075\lambda$         &  $0.05\lambda$          & 
$0.025\lambda$             \\
\hline
$e$   &  $8.38\times 10^{-1}$   &  $8.06\times 10^{-1}$   &  $5.56\times
10^{-1}$   &  $2.95\times 10^{-1}$   \\
\hline
\hline
\end{tabular}}
\end{center}
\end{table}

{\em Example 2}. This example illustrates the reconstruction results of a
non-smooth grating profile with infinitely many Fourier modes, as seen in Fig.
\ref{fig:pro}(b). The exact grating surface function is given by
$\phi(\rho)=\de\ps(\rho)$, where the grating profile function
\[
\ps(x,y)=|\cos(2\pi x)\cos(2\pi y)| - |\sin(\pi x)\sin(2\pi y)|.
\]
First is to consider the influence of $\de$ by using noise-free data. The
measurement is taken at $h = 0.1\lambda$. Table \ref{tab:ex2} presents the
relative $L^2(R)$ error of the reconstructions with three different
surface deformation parameters $\de = 0.05\lambda, 0.025\lambda, 0.0125\lambda$.
The error decreases from $35.8\%$ to $16.0\%$ as $\de$ decreases from
$0.05\lambda$ to $0.0125\lambda$. Based on these results, the following
observation can be made: a smaller deformation parameter $\de$ yields a better
reconstruction.

\begin{table}[http!]
\caption{Example 2: Relative error of the reconstructions by using different
$\de$ with $h=0.1\lambda$
and $\ga = 0.0$.}\label{tab:ex2}
\begin{center}{
\begin{tabular}{llll}
\hline
\hline
$\de$  &  $0.05\lambda$         &  $0.025\lambda$         &  $0.0125\lambda$    
    \\
\hline
$e$    &  $3.58\times 10^{-1}$  &  $2.72\times 10^{-1}$  &  $1.60\times 10^{-1}$
     \\
\hline
\hline
\end{tabular}}
\end{center}
\end{table}

\begin{figure}
\center
\includegraphics[width=0.4\textwidth]{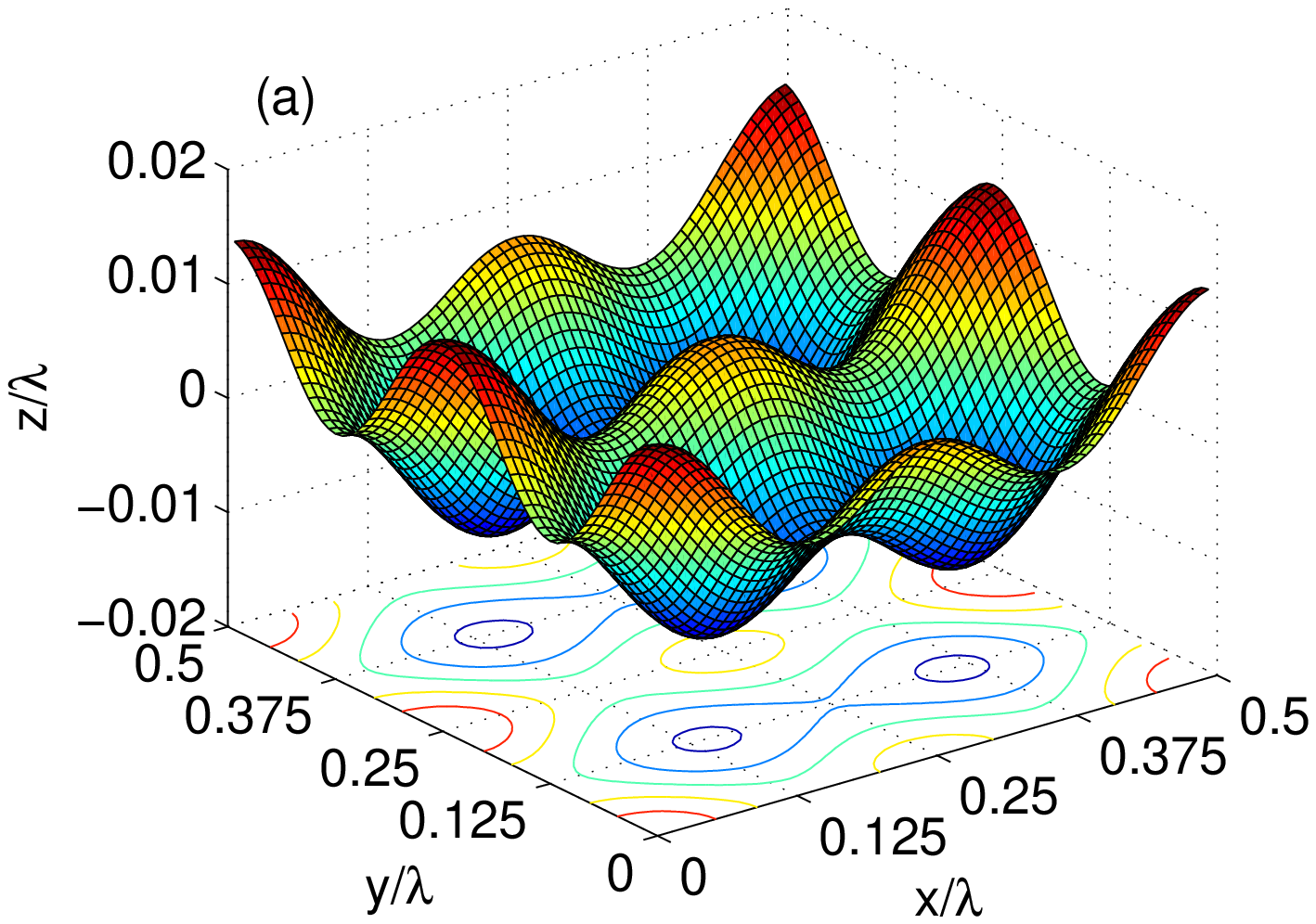}
\includegraphics[width=0.4\textwidth]{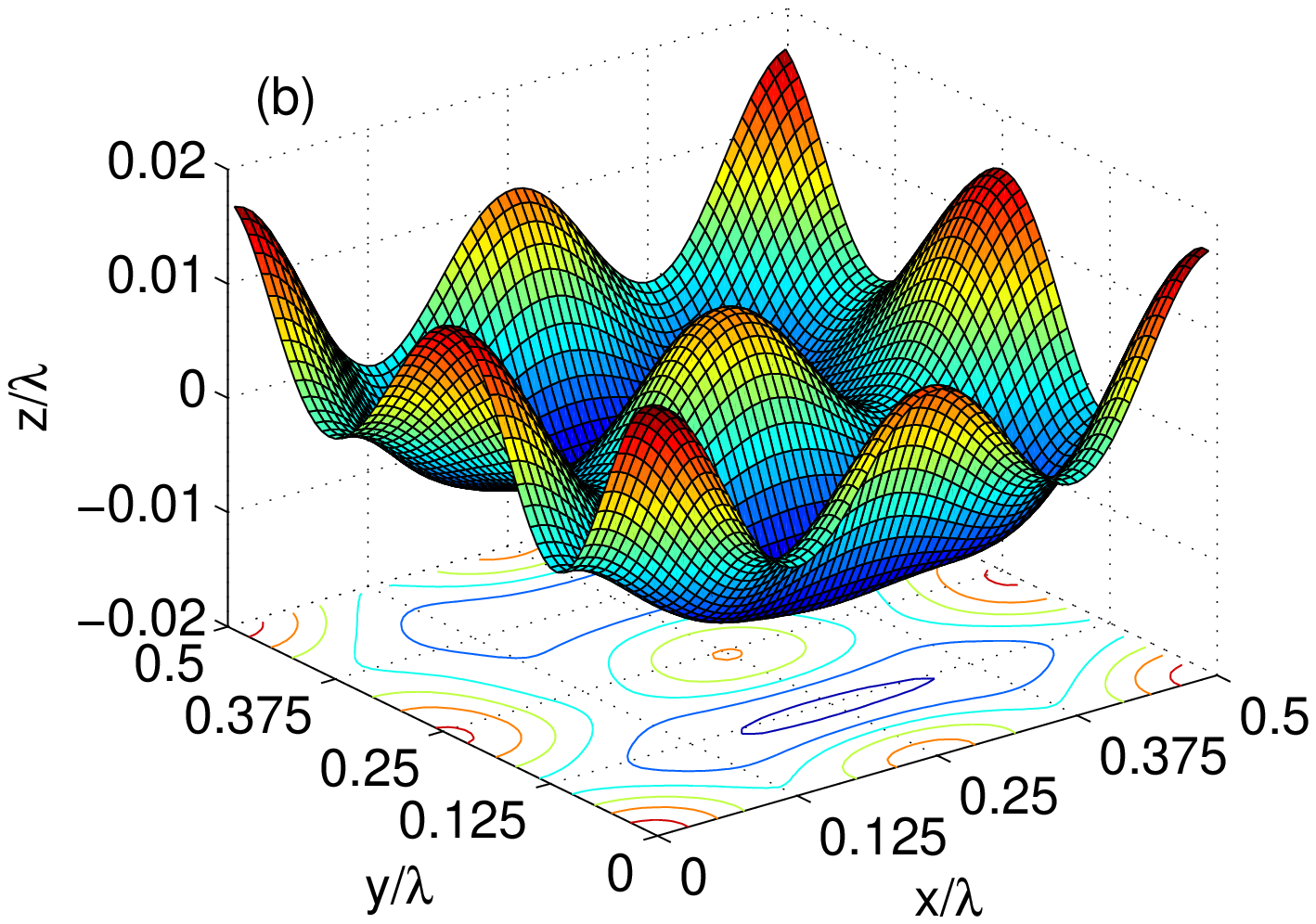}
\includegraphics[width=0.4\textwidth]{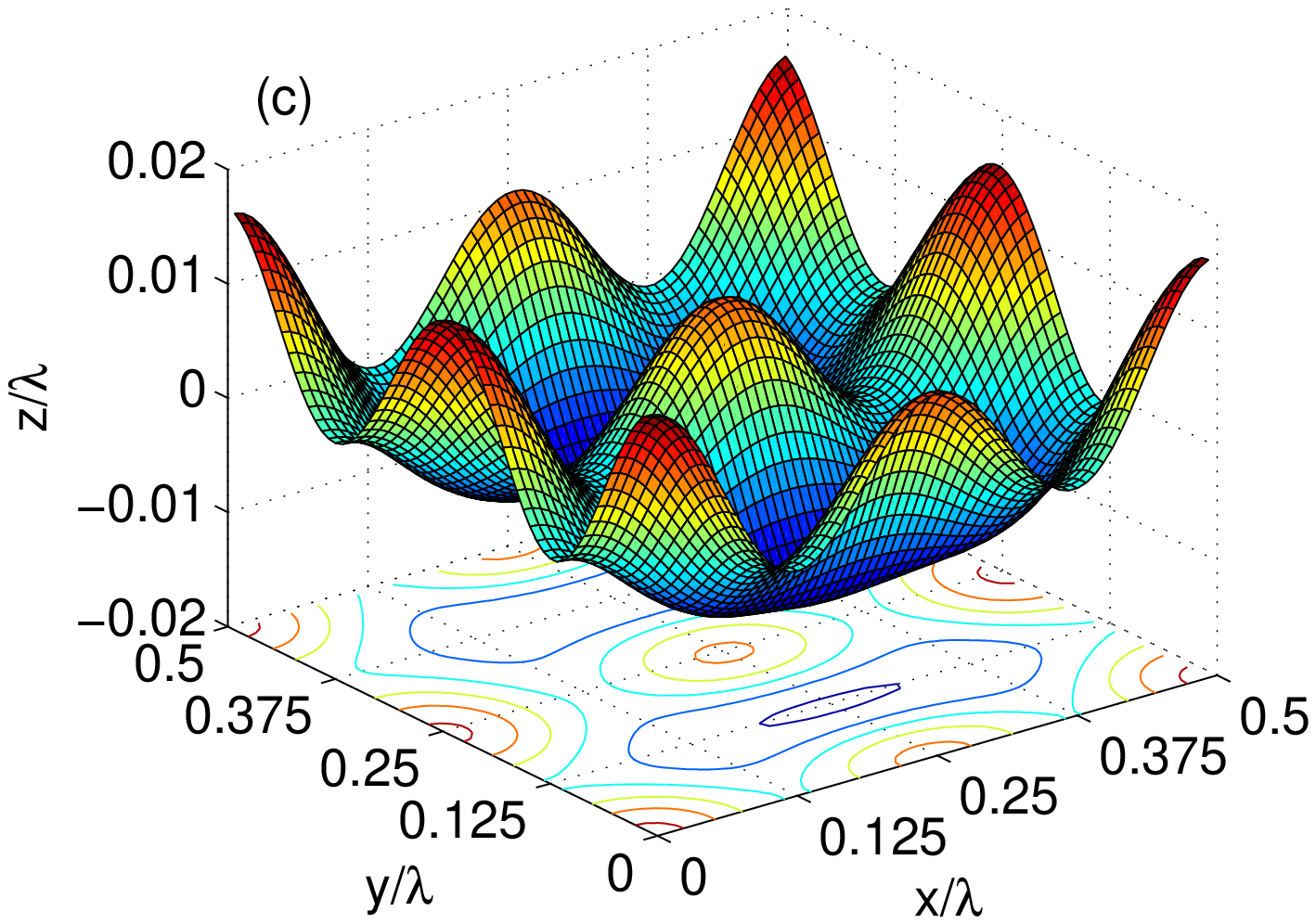}
\includegraphics[width=0.4\textwidth]{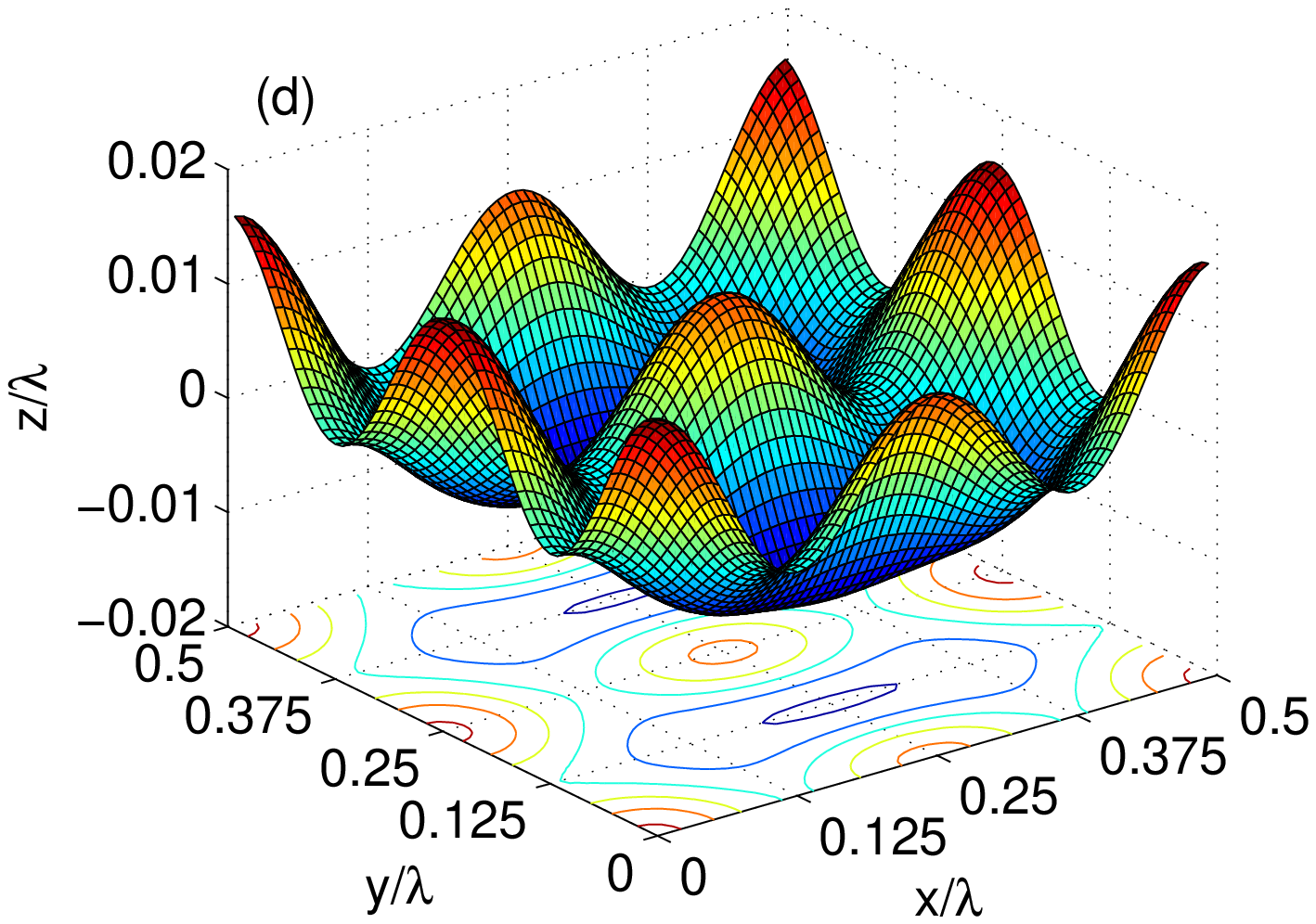}
\caption{Example 2: Reconstructed grating surfaces by using different $h$ with
$\delta = 0.0125\lambda$ and $\ga = 1\%$. (a) $h = 0.1\lambda$; (b) $h =
0.075\lambda$; (c) $h = 0.05\lambda$; (d) $h = 0.025\lambda$.}
\label{fig:ex2}
\end{figure}

Next is to consider the influence of the noise level $\ga$ and the measurement
distance $h$. We add $1\%$ and $5\%$ random noises to the scattering data.
Table \ref{tab:ex2-1} and \ref{tab:ex2-2} report
the relative $L^2(R)$ error of the reconstructions with four different
measurement distances $h =0.1\lambda, 0.075\lambda, 0.05\lambda, 0.025\lambda$
for a fixed $\de = 0.0125\lambda$. Comparing the results for the same $\de =
0.0125\lambda$ and $h =0.1\lambda$ in Tables \ref{tab:ex2} and \ref{tab:ex2-2},
we can see that the relative error is more than doubled from $16.0\%$ by using
noise-free data to $34.3\%$ by using $5\%$ noise data. Again, the reason is that
a smaller cut-off is chosen to suppress the exponentially increasing noise in
the data and thus higher Fourier modes of the exact grating surface function can
not be recovered. A smaller measurement distance helps to enhance the resolution
and reduce the error. In Table \ref{tab:ex2-1}, the reconstruction error
decreases from $27.3\%$ by using $h = 0.1\lambda$ to as low as $17.3\%$ by using
$h=0.025\lambda$ for $1\%$ noise data. In Table \ref{tab:ex2-2}, the
reconstruction error decreases from $34.3\%$ by using $h = 0.1\lambda$ to as low
as $24.4\%$ by using $h=0.025\lambda$ for $5\%$ noise data. Figure \ref{fig:ex2}
shows the reconstructed surfaces by using $h =0.1\lambda, 0.075\lambda,
0.05\lambda, 0.025\lambda$. Comparing the exact surface profile
in Fig. \ref{fig:pro}(b) and the reconstructed surface in Fig. \ref{fig:ex2}(d),
we can see that a good reconstruction can still be obtained when using a small
measurement distance.

\begin{table}[http!]
\caption{Example 2: Relative error of the reconstructions by using different $h$
with $\de=0.0125\lambda$ and $\ga = 1\%$.}\label{tab:ex2-1}
\begin{center}{
\begin{tabular}{lllll}
\hline
\hline
$h$     &  $0.1\lambda$       &  $0.075\lambda$         &  $0.05\lambda$        
&  $0.025\lambda$     \\
\hline
$e$    & $2.73\times 10^{-1}$ & $2.44\times 10^{-1}$    &  $1.88\times 10^{-1}$ 
&  $1.73\times 10^{-1}$   \\
\hline
\hline
\end{tabular}}
\end{center}
\end{table}

\begin{table}[http!]
\caption{Example 2: Relative error of the reconstructions by using different $h$
with $\de=0.0125\lambda$ and $\ga = 5\%$.}\label{tab:ex2-2}
\begin{center}{
\begin{tabular}{lllll}
\hline
\hline
$h$    &  $0.1\lambda$     &  $0.075\lambda$           &  $0.05\lambda$         
&  $0.025\lambda$             \\
\hline
$e$    &  $3.43\times 10^{-1}$ &  $2.99\times 10^{-1}$     &  $2.81\times
10^{-1}$   &  $2.44\times 10^{-1}$   \\
\hline
\hline
\end{tabular}}
\end{center}
\end{table}

\section{Conclusion}

We have presented an effective computational method to reconstruct surfaces of
biperiodic dielectric gratings. Subwavelength resolution is achieved stably.
Based on the transformed field expansion, an analytic solution is deduced for
the direct problem. The nonlinear inverse problem is linearized by dropping
higher order terms in power series. Explicit reconstruction formulas are
obtained and are implemented by using the FFT. Two representative numerical
examples are considered: one is a smooth function which has finitely many
Fourier modes and another is a nonsmooth function which has infinitely many
Fourier modes. We have carefully investigated the influence of the
parameters on the reconstructions. The results show that super resolution may be
achieved by using small measurement distance. There are many interesting and
challenging mathematical problems, such as uniqueness, stability, resolution,
and error estimates, which are remaining and left for future work. We will
report the results elsewhere.


\begin{thebibliography}{10}

\bibitem{mumps1}
P. R. Amestoy, I. S. Duff, J. Koster, and J.-Y. L'Excellent, {\em A
fully asynchronous multifrontal solver using distributed dynamic
scheduling}, SIAM. J. Matrix Anal. \& Appl., 23 (2001), pp. 15--41.

%\bibitem{mumps2}
%P. R. Amestoy and A. Guermouche and J.-Y. L'Excellent, and S. Pralet, {\em
%Hybrid scheduling for the parallel solution of linear systems}, Parallel
%Computing, 32 (2006), pp. 136--156.

\bibitem{Amm}
{\sc H. Ammari}, {\em Uniqueness theorems for an inverse problem in a doubly
periodic structure}, Inverse Problems, 11 (1995), pp. 823--833.


%\bibitem{A-G-S-1}
%{\sc H. Ammari, J. Garnier, and K. Solna}, {\em Resolution and stability
%analysis in full-aperture, linearized conductivity and wave imaging}, Proc.
%Amer. Math. Soc., 141 (2013), pp. 3431--3446.

%\bibitem{A-G-S-2}
%{\sc H. Ammari, J. Garnier, and K. Solna}, {\em Partial data resolving power of
%conductivity imaging from boundary measurements}, SIAM J. Math. Anal., 45
%(2013), pp. 1704--1722.


\bibitem{A-K-Y}
{\sc I. Akduman, R. Kress, and A. Yapar}, {\em Iterative reconstruction of
dielectric rough surface profiles at fixed frequency}, Inverse Problems, 22
(2006), pp. 939--954.

\bibitem{A-K}
{\sc T. Arens and A. Kirsch}, {\em The factorization method in inverse
scattering from periodic structures}, Inverse Problems, 19 (2003), pp.
1195--1211.

\bibitem{Bao-1}
{\sc G. Bao}, {\em A unique theorem for an inverse problem in periodic
diffractive optics}, Inverse Problems, 10 (1994), pp. 335--340.

\bibitem{Bao-2}
{\sc G. Bao}, {\em Variational approximation of Maxwell's equations in
biperiodic structures}, SIAM J. Appl. Math., 57 (1997), pp. 364--381.

\bibitem{B-C-M}
{\sc G. Bao, L. Cowsar, and W. Masters}, {\em Mathematical Modeling in Optical
Science}, Frontiers Appl. Math., 22, SIAM, Philadelphia, 2001.

\bibitem{B-C-L}
{\sc G. Bao, T. Cui, and P. Li}, {\em Inverse diffraction grating of Maxwell's
equations in biperiodic structures}, Optics Express, 22 (2014), pp. 4799--4816.

\bibitem{B-D-C}
{\sc G. Bao, D. Dobson, and J. A. Cox}, {\em Mathematical studies in rigorous
grating theory}, J. Opt. Soc. Amer. A, 12 (1995), pp. 1029--1042 .

\bibitem{B-F}
{\sc G. Bao and A. Friedman}, {\em Inverse problems for scattering by periodic
structure}, Arch. Rational Mech. Anal., 132 (1995), pp. 49--72.

\bibitem{B-Li-1}
{\sc G. Bao and P. Li}, {\em Near-field imaging of infinite rough surfaces},
SIAM J. Appl. Math., 73 (2013), pp. 2162--2187.

\bibitem{B-Li-2}
{\sc G. Bao and P. Li}, {\em Near-field imaging of infinite rough surfaces in
dielectric media}, SIAM J. Imaging Sci.,  7 (2014), pp. 867--899.

\bibitem{B-Li-3}
{\sc G. Bao and P. Li}, {\em Convergence analysis in near-field imaging},
Inverse Problems, 30 (2014), pp. 085008.  

\bibitem{B-Lin}
{\sc G. Bao and J. Lin}, {\em Near-field imaging of the surface displacement on
an infinite ground plane}, Inverse Probl. Imag., 7 (2013), pp. 377--396.

\bibitem{B-L-W-1}
{\sc G. Bao, P. Li, and H. Wu}, {\em An adaptive edge element method with
perfectly matched absorbing layers for wave scattering by biperiodic
structures}, Math. Comp., 79 (2009), pp. 1--34.

\bibitem{B-L-W-2}
{\sc G. Bao, P. Li, and H. Wu}, {\em A computational inverse diffraction grating
problem}, J. Opt. Soc. Am. A, 29 (2012), pp. 394--399.

\bibitem{B-L-L}
{\sc G. Bao, P. Li, and J. Lv}, {\em Numerical solution of an inverse
diffraction grating problem from phasless data}, J. Opt. Soc. Am. A, 30 (2013),
pp. 293--299.

\bibitem{B-Z-Z}
{\sc G. Bao, H. Zhang, and J. Zou}, {\em Unique determination of periodic
polyhedral structures by scattered electromagnetic fields}, Trans. Amer. Math.
Soc., 363 (2011), pp. 4527--4551.

\bibitem{B-Z}
{\sc G. Bao and Z. Zhou}, {\em An inverse problem for scattering by a doubly
periodic structure}, Trans. Amer. Math. Soc., 350 (1998), 4089--4103.

\bibitem{B-C-Y}
{\sc G. Bruckner, J. Cheng, and M. Yamamoto}, {\em An inverse problem in
diffractive optics: conditional stability}, Inverse Problems, 18 (2002), pp.
415--433.

\bibitem{B-E}
{\sc G. Bruckner and J. Elschner}, {\em A two-step algorithm for
the reconstruction of perfectly reflecting periodic profiles}, Inverse Problems,
19 (2003), pp. 315--329.

\bibitem{B-R-1}
{\sc O. Bruno and F. Reitich}, {\em Numerical solution of diffraction problems:
a method of variation of boundaries}, J. Opt. Soc. Am. A, 10 (1993), pp.
1168-1175.

%\bibitem{B-R-2}
%{\sc O. Bruno and F. Reitich}, {\em Numerical solution of diffraction problems:
%a method of variation of boundaries. III. Doubly periodic gratings}, J. Opt.
%Soc. Am. A, 10 (1993), pp. 2551--2562.

\bibitem{C-S-1}
{\sc S. Carney and J. Schotland}, {\em Inverse scattering for near-field
microscopy}, App. Phys. Lett., 77 (2000), pp. 2798--2800.

\bibitem{C-S-2}
{\sc S. Carney and J. Schotland}, {\em Near-field tomography}, MSRI Ser. Math.
Appl., 47 (2003), pp. 133--168.

%\bibitem{C-W}
%{\sc Z. Chen and H. Wu}, {\em An adaptive finite element method with
%perfectly matched absorbing layers for the wave scattering by periodic
%structures}, SIAM J. Numer. Anal., 41 (2003), pp. 799--826.

\bibitem{C-L-W}
T. Cheng, P. Li, and Y. Wang, {\em Near-field imaging of perfectly conducting
grating surfaces}, J. Opt. Soc. Am. A, 30 (2013), pp. 2473--2481.

\bibitem{C-G-H-I-R}
{\sc R. Coifman, M. Goldberg, T. Hrycak, M. Israeli, and V. Rokhlin}, {\em An
improved operator expansion algorithm for direct and inverse scattering
computations}, Waves Random Media, 9 (1999), pp. 441--457.

\bibitem{Cou}
{\sc D. Courjon}, {\em Near-Field Microscopy and Near-Field Optics}, Imperial
College Press, London, 2003.

\bibitem{D-W}
{\sc J.~A. DeSanto and R.~J. Wombell}, {\em The reconstruction of shallow
rough-surface profiles from scattered field data}, Inverse Problems, 7 (1991),
pp. L7--L12.

%\bibitem{Dob-1}
%{\sc D. Dobson}, {\em Optimal design of periodic antireflective structures for
%the Helmholtz equation}, Eur. J. Appl. Math., 4 (1993), pp. 321--340.

%\bibitem{Dob-2}
%{\sc D. Dobson}, {\em Optimal shape design of blazed diffraction grating}, J.
%Appl. Math. Optim., 40 (1999), pp. 61--78.

\bibitem{Dob-3}
{\sc D. Dobson}, {\em A variational method for electromagnetic diffraction in
biperiodic structures}, Math. Model. Numer. Anal., 28 (1994), pp. 419--439.

\bibitem{E-H-R}
{\sc J. Elschner, G. Hsiao, and A. Rathsfeld}, {\em Grating
profile reconstruction based on finite elements and optimization techniques},
SIAM J. Appl. Math., 64 (2003), pp. 525--545.

%\bibitem{E-S-1}
%{\sc J. Elschner and G. Schmidt}, {\em Diffraction in periodic structures and
%optimal design of binary gratings: I. Direct problems and gradient formulas},
%Math. Methods. Appl. Sci., 21 (1998), pp. 1297--1342.

%\bibitem{E-S-2}
%{\sc J. Elschner and G. Schmidt}, {\em Numerical solution of optimal design
%problems for binary gratings}, J. Comput. Phys., 146 (1998), pp. 603--626.

%\bibitem{G-N}
%{\sc N. Garc\'{i}a and M. Nieto-Vesperinas}, {\em Near-field
%optics inverse-scattering reconstruction of reflective surfaces}, Optics
%Letters, 24 (1993), pp. 2090--2092.

%\bibitem{H-N-S}
%{\sc Y. He, D.~P. Nicholls, and J. Shen}, {\em An efficient and stable spectral
%method for electromagnetic scattering from a layered periodic struture}, J.
%Comput. Phys., 231 (2012), pp. 3007--3022.

\bibitem{Het}
{\sc F. Hettlich}, {\em Iterative regularization schemes in inverse scattering
by periodic structures}, Inverse Problems, 18 (2002), pp. 701--714.

\bibitem{H-K}
{\sc F. Hettlich and A. Kirsch}, {\em Schiffer's theorem in inverse scattering
theory for periodic structures}, Inverse Problems, 13 (1997), pp. 351--361.

\bibitem{H-Y-Z}
{\sc G. Hu, J. Yang, and B. Zhang}, {\em An inverse electromagnetic scattering
problem for a bi-periodic inhomogeneous layer on a perfectly conducting plate},
Appl. Anal., 90 (2011), 317--333.

\bibitem{H-Z}
{\sc G. Hu and B. Zhang}, {\em The linear sampling method for
inverse electromagnetic scattering by a partially coated bi-periodic
structures}, Math. Meth. Appl. Sci., 34 (2011), pp. 509--519.

\bibitem{I-R}
{\sc K. Ito and F. Reitich}, {\em A high-order perturbation approach to profile
reconstruction: I. Perfectly conducting gratings}, Inverse Problems, 15 (1999),
pp. 1067--1085.

\bibitem{Kir}
{\sc A. Kirsch}, {\em Uniqueness theorems in inverse scattering theory for
periodic structures}, Inverse Problems, 10 (1994), pp. 145--152.

\bibitem{K-T}
{\sc R. Kress and T. Tran}, {\em Inverse scattering for a locally perturbed
half-plane}, Inverse Problems, 16 (2000), pp. 1541--1559.

\bibitem{L-N-1}
{\sc A. Lechleiter and D.~L. Nguyen}, {\em On uniqueness in electromagnetic
scattering from biperiodic structures}, ESAIM: M2AN, 47 (2013), pp. 1167--1184.

\bibitem{L-N-2}
{\sc A. Lechleiter and D.~L. Nguyen}, {\em Factorization method for
electromagnetic inverse scattering from biperiodic structures}, SIAM J. Imaging
Sci., 6 (2013), pp. 1111-1139.

%\bibitem{M-N-1}
%{\sc A. Malcolm and D.~P. Nicholls}, {\em A field expansions method for
%scattering by periodic multilayered media}, J. Acout. Soc. Am., 129 (2011), pp.
%1783--1793.

%\bibitem{M-N-2}
%{\sc A. Malcolm and D.~P. Nicholls}, {\em A boundary perturbation method for
%recovering interface shapes in layered media}, Inverse Problems, 27 (2011), pp.
%095009.

%\bibitem{M-N-3}
%{\sc A. Malcolm and D.~P. Nicholls}, {\em Operator expansions and constrained
%quadratic optimization for interface reconstruction: Impenetrable periodic
%acoustic media}, Wave Motion, to appear.

\bibitem{N-S}
{\sc J. C. N\'{e}d\'{e}lec and F. Starling}, {\em Integral equation methods in a
quasi-periodic diffraction problem for the time-harmonic Maxwell's equations},
SIAM J. Math. Anal., 22 (1991), pp. 1679--1701.

%\bibitem{Ngu}
%{\sc D.~L. Nguyen}, {\em Spectral Methods for Direct and Inverse Scattering
%from Periodic Structures}, PhD thesis, Ecole Polytechnique, Palaiseau, France,
%2012.

\bibitem{N-R-1}
{\sc D.~P. Nicholls and F. Reitich}, {\em Shape deformations in rough surface
scattering: cancellations, conditioning, and convergence}, J. Opt. Soc. Am. A,
21 (2004), pp. 590-605.

%\bibitem{N-R-2}
%{\sc D.~P. Nicholls and F. Reitich}, {\em Shape deformations in rough surface
%scattering: improved algorithms}, J. Opt. Soc. Am. A, 21 (2004), pp. 606-621.

\bibitem{Pet}
{\sc R. Petit}, ed., {\em Electromagnetic Theory of Gratings}, Springer-Verlag,
1980.

\bibitem{phg}
{\em PHG (Parallel Hierarchical Grid)}, \texttt{http://lsec.cc.ac.cn/phg/}.

%\bibitem{San}
%K. Sandfort, {\em The factorization method for inverse scattering from periodic
%inhomogeneous media}, PhD thesis, Karlsruher Institut f\"{u}r Technologie,
%2010.

%\bibitem{W-L}
%{\sc Y. Wu and Y.~Y.  Lu}, {\em Analyzing diffraction gratings by a boundary
%integral equation Neumann-to-Dirichlet map method}, J. Opt. Soc. Am. A, 26
%(2009), pp. 2444--2451.

\bibitem{Y-Z}
{\sc J. Yang and B. Zhang}, {\em Inverse electromagnetic scattering problems by
a doubly periodic structure}, Math. Appl. Anal., 18 (2011), pp. 111--126.

\end{thebibliography}
\end{document}